\documentclass{article}
\usepackage[pdftex]{graphicx}
\usepackage{amsmath, amsthm, amssymb, url, enumerate}
\usepackage[all]{xy}
\usepackage{cases}
\usepackage{stmaryrd}

\newtheorem{Def}{Definition}[section]
\newtheorem{Thm}[Def]{Theorem}
\newtheorem{aaa}{Proposition 3.1}

\newtheorem{bbb}{Proposition 3.2}

\newtheorem{ccc}{Proposition 3.3}

\newtheorem{Prop}[Def]{Proposition}

\newtheorem{Rmk}[Def]{Remark}

\newtheorem{Ex}[Def]{Example}

\makeatletter
\def\mapstofill@{%
   \arrowfill@{\mapstochar\relbar}\relbar\rightarrow}
\newcommand*\xmapsto[2][]{%
   \ext@arrow 0395\mapstofill@{#1}{#2}}
\makeatother

\title{Crepant resolutions and ${\rm Hilb}^G(\mathbb{C}^4)$ for certain abelian subgroups of $SL(4,\mathbb{C})$}

\author{Yusuke Sato}
\date{}
\begin{document}
\maketitle

\begin{center}
\textbf{Abstract.}
\end{center}
Let $G$ be a finite subgroup of $SL(n,\mathbb{C})$, then the quotient $\mathbb{C}^n/G$ has a Gorenstein canonical singularity. In this paper, we will show several examples of crepant resolutions in dimension $4$ and show examples in which ${\rm Hilb}^G(\mathbb{C}^4)$ is blow-up of certain crepant resolutions for $\mathbb{C}^4/G$, or ${\rm Hilb}^G(\mathbb{C}^4)$ has singularity.
\vspace{1cm}

\section{Introduction}

Let $G$ be a finite subgroup of $SL(n,\mathbb{C})$, then the quotient $\mathbb{C}^n/G$ has a Gorenstein canonical singularity. In paticular, if $n=2$,  then $\mathbb{C}^2/G$ is hypersurface and has singularity which is called a rational double point or ADE singularity. 
A {\it crepant resolutions} is defined as follows.
\begin{Def}
Let $f:Y \to X$  be a resolution, then the adjunction formula  $K_Y=f^*K_X+\sum_{i=1}^n a_iD_i$ has a rational number $a_i$.\\
 If $f$ satisfies $a_i=0$ for all $i$, then $f$ is called a crepant resolution. 
\end{Def}

It is well known that $\mathbb{C}^2/G$ and
$\mathbb{C}^3/G$ have a crepant resolutions(see [7]). However, crepant resolutions does not necessarily exist in dimension $\geq 4$. Dais-Henk-Ziegler([3]) shows some conditions on $G$ to have a crepant resolution of the abelian groups of $SL(n,\mathbb{C})$, and Hayashi-Ito-Sekiya([5]) show some examples of non abelian subgruops of $SL(4,\mathbb{C})$. The criterion to have a crepant resolution has not been found yet.\\
Ito and Nakamura introduced the $G$-Hilbert scheme for finite groups $G$ of $SL(2,\mathbb{C})$, and they proved that $G$-Hilbert scheme is the minimal resolution(that is crepant resolution) of $\mathbb{C}^2/G$ ([8]). In dimension three, crepant resolution is not unique, and $G$-Hilbert scheme is a one of crepant resolutions([1]).
However, in the case $n\geq4$, even if assuming $\mathbb{C}^4/G$ has a crepant resolution, the relationship between ${\rm Hilb}^G(\mathbb{C}^n)$ and $\mathbb{C}^n/G$ is not well known. Table $1$ shows previous stidies about ${\rm Hilb}^G(\mathbb{C}^n)$ and $\mathbb{C}^n/G$\\
\small{
\begin{table}[htb]
\begin{center}
\caption{The relationship $\mathbb{C}^n/G$ and ${\rm Hilb}^G(\mathbb{C}^n)$}
 \begin{tabular}{|l||c|c|c|} \hline
    　　 & $SL(n,\mathbb{C})$              &  $GL(n,\mathbb{C})$\\ \hline \hline
    n=2 & crepant resolution(Ito-Nakamura[8]) &  minimal resolution([6],[10])\\ \hline
    n=3 & one of crepant resolutions             &    \\
    &     (Nakamura[12], Bridgeland-King-Reid[1]) &   There is example  \\ \cline{1-2} 
    n=4 & There is examples  in which    & in which\\
     &  ${\rm Hilb}^G(\mathbb{C}^n)$ is  not crepant resolutions & ${\rm Hilb}^G(\mathbb{C}^n)$ is singular \\
     & For example :$G=\frac{1}{2}(1,1,1,1)$ & \\ \hline
  \end{tabular}\\
\end{center}
\end{table} 
 }
 In this paper, we will introduce two series of abelian subgroups which has crepant resolutions but ${\rm Hilb}^G(\mathbb{C}^4)$ is not a crepant resolution. 
 
\begin{aaa}\label{m1}
If $G$ is the following type, then $\mathbb{C}^4/G$ has crepant resolutions.
\begin{itemize}
\item[{\rm(i)}] $G=\left< \frac{1}{r}(1,1,0,r-2),\frac{1}{r}(0,0,1,r-1)\right>$.
\item[{\rm(ii)}] $G=\frac{1}{r}(1,a,a^2,a^3)$, where $r=1+a+a^2+a^3$
\end{itemize}
\end{aaa}

\begin{bbb}\label{m2}
If $G$ is the type (i) or (ii), the relationship between ${\rm Hilb}^G(\mathbb{C}^n)$ and $\mathbb{C}^n/G$ is as follows.
\begin{itemize}
\item[{\rm(i)}] If $G=\left< \frac{1}{r}(1,1,0,r-2),\frac{1}{r}(0,0,1,r-1)\right>$  and $r$ is odd, then ${\rm Hilb}^G(\mathbb{C}^4)$ is blow-up of certain crepant resolutions. If $r$ is even, then ${\rm Hilb}^G(\mathbb{C}^4)$ is one of crepant resolutions. 
\item[{\rm(ii)}] If $G=\frac{1}{r}(1,a,a^2,a^3)$ and  if $a=3$, then ${\rm Hilb}^G(\mathbb{C}^4)$ is blow-up of certain crepant resolutions.
\end{itemize}
\end{bbb}

In addition to the above example,  we will show a example in which ${\rm Hilb}^G(\mathbb{C}^4)$ is singular.In this case, $\mathbb{C}^n/G$ has not crepant resolutions.

\begin{ccc}\label{m3}
If $G$ is generated by $g=\frac{1}{2m}(1,2m-1,m,m)$, then ${\rm Hilb}^G(\mathbb{C}^4)$ has singularity.  
\end{ccc}

\section{Notation}
In this section we set up the notation of toric geometry and definition of ${\rm Hilb}^G(\mathbb{C}^n)$. We use the same notation as \cite{5}.
\subsection{Notation}
In this subsection, we introduce quotient singularity $\mathbb{C}^n/G$ and crepant resluitons as toric varieties.
Assume that $G$ is a finite abelian subgroup of $SL(n,\mathbb{C})$.
Let order of $G$ is $r$, then any $g \in G$ is of the form $g ={\rm diag}(\varepsilon^{a_1}_r, \dots,\varepsilon^{a_n}_r)$, where $\varepsilon_r$ is a primitive $r$th root of unity. Then we can write $g=\frac{1}{r}(a_1,\dots,a_n)$.
In addition, the map $\phi:G \to \mathbb{R}^n$ is defined by $\phi(g)=\bar{g}=\frac{1}{r}(a_1,\dots,a_n) $ for $g \in G$.\\
Let  $N :=\mathbb{Z}^n+\sum_{g \in G} \mathbb{Z}\bar{g}$ be a free $\mathbb{Z}$-module of rank $n$, $M$ be the dual $\mathbb{Z}$-module of $N$, and $N_{\mathbb{R}}=N \otimes_{\mathbb{Z}} \mathbb{R}$, $M_{\mathbb{R}}=M \otimes_{\mathbb{Z}} \mathbb{R}$.\\
We will denote by $\sigma$ the region of $\mathbb{R}^n$ whose all entries are non-negative.\\
Then the toric variety $U_{\sigma}:={\rm Spec}\mathbb{C}[\sigma^{\vee}\cap M]$ is isomorphic to $\mathbb{C}^n/G$
We recall some of important definition for toric quotient singularity as follows.

\begin{Def}\label{smoothness}
{\rm The cone $\sigma$ is} smooth {\rm if its generators is a part of basis of $N$. Also, the fan $\Delta$ is} smooth {\rm if its all cones are smooth.}
\end{Def}

\begin{Def}\label{age}
{\rm Let $\bar{g}=\frac{1}{r}(a_1,\dots,a_n)$ be a lattice point of $N \cap \sigma$. The} age {\rm of $\bar{g}$ is defined by age($\bar{g}$)$=\frac{1}{r}\sum_{i=1}^n a_i$. Since $G$ is subgruop of $SL(n,\mathbb{C})$,} age {\rm is  an positive integer. }
\end{Def}

\begin{Rmk}
{\rm Let $\Delta$ is subdivision of $\sigma$ using by lattice points of age($\bar{g}$)=1. If the toric variety $Y$ determined by  $\Delta$ is smooth, then $Y$ is a crepant resolution of  $\mathbb{C}^n/G$}
\end{Rmk}

In the case of dimension $4$, the lattice points of age one are on the tetrahedron with verticies $e_1=(1,0,0,0),\dots, e_4=(0,0,0,1)$. So we consider subdivison on the tetrahedron. Also, we call this tetrahedron junior simplex.

\subsection{$G$-graph and ${\rm Hilb}^G(\mathbb{C}^n)$} 
The purpose of this subsection is to recall a way of construction ${\rm Hilb}^G(\mathbb{C}^n) $ via $G$-graph.\\
First, we introduce definition of ${\rm Hilb}^G(\mathbb{C}^n)$.
Let ${\rm Hilb}^r(\mathbb{C}^n)$ is Hilbert scheme of $r$-points.
Assume that order of $G$ equals to $r$, then $G$ acts on ${\rm Hilb}^r(\mathbb{C}^n)$ and $S^r(\mathbb{C}^n)$.
${\rm Hilb}^G(\mathbb{C}^n)$ is the irreducible component of $({\rm Hilb}^r(\mathbb{C}^n)^G$ which dominates $\mathbb{C}^n/G$, and ${\rm Hilb}^G(\mathbb{C}^n)$ is birational to $\mathbb{C}^n/G$ via the Hilbert-Chow morphism.
$$\begin{matrix}
{\rm Hilb}^r(\mathbb{C}^n) & \xrightarrow{{\rm H-C\  morphism}} & S^r(\mathbb{C}^n) \\
  \cup & & \cup \\
{\rm Hilb}^r(\mathbb{C}^n)^G & \xrightarrow{\qquad \quad \qquad} & S^r(\mathbb{C}^n)^G  \\
\cup & & \cup \\
{\rm Hilb}^G(\mathbb{C}^n)  &\xrightarrow{\quad {\rm dominate} \quad}& \mathbb{C}^n/G \\
\end{matrix}$$
Second, we recall definition of $G$-graph
Let $S=\mathbb{C}[x_1,\dots,x_n]$  denote the coordinate ring of $\mathbb{C}^n$ and $\mathcal{M}$ be the set of all monomials in $S$ and $1$, and $\rho_i$ be irreducible representation of $G$. We will denote by $X^u$ a monomial in $\mathcal{M}$ where $X^u=x_1^{u_1}\cdots x_n^{u_n}$ and $u=(u_1, \dots, u_n) \in \mathbb{Z}_{\geq0}^n$.\\

We write ${\rm wt}(X^u)=\rho_i$ if $X^u(g\cdot p)=\rho_i(g)X^u(p)$ holds for any $g \in G$ and $p \in \mathbb{C}^n$. Since any monomial is contained in some $\rho_i$, we can define a map ${\rm wt}: \mathcal{M} \to {\rm Irr}(G)$, where ${\rm Irr}(G)$ os the set of irreducible representation of $G$.
In this paper, we define $G$-graph using by the map wt and ideal of $S$.
\begin{Def}\label{ggraph}
{\rm Let $I \subset S$ be an ideal, we define a subset $\Gamma(I) \subset \mathcal{M}$ such that $\{X \in \mathcal{M} \mid X \notin I \} $.
A Subset $\Gamma(I)$ is called a} $G$-graph {\rm if the restriction map ${\rm wt}: \Gamma(I) \to {\rm Irr}(G)$ is a bijection.}
\end{Def}

We mention that the condition $\{X \in \mathcal{M} \mid X \notin I \} $ is equivalent to $X \in \Gamma$ and $X$ is divided by $Y \in \mathcal{M}$, then $Y \in \Gamma$ .

\begin{Def}\label{sigma}
{\rm Let $A_{\Gamma}$ be a set of minimal generaters of $I(\Gamma)$. 
We difine the map ${\rm wt}_{\Gamma}:\mathcal{M} \to G-{\rm graph}$ as ${\rm wt}_{\Gamma}(X^u)=\bar{X^u}$ such that ${\rm wt}(X^u)={\rm wt}(\bar{X^u})$ .
For a $G$-graph $\Gamma$, we define the rational cone
$$
\sigma(\Gamma):=\{w \in N_{\mathbb{R}} \mid w \cdot X^u > w \cdot {\rm wt}_{\Gamma}(X^u) \  for\  all\  X^u \in A_{\Gamma}\},
$$
where $w \cdot X^u$ means standard inner product $w \cdot u$ in $\mathbb{R}^n$.  }
\end{Def}

We will denote by $Fan(G)$ the fan in $N_{\mathbb{R}}$ is defined by all closed cone $\sigma(\Gamma)$ and all their faces. Followong theorem says that we can calculate ${\rm Hilb}^G$ using by $G$-graph.

\begin{Thm}\label{Nakamura}{\rm ([Nakamura], Theorem 2,11) The following hold.
\begin{itemize}
\item[(i)] Fan($G$) is a finite fan with its support $\Delta$
\item[(ii)] The normalization of ${\rm Hilb}^G$ is isomorphic to $T(Fan(G))$
\end{itemize}
}
\end{Thm}

We show twoexamples of the ${\rm Hilb}^G(\mathbb{C}^n)$ in dimension two or three.
\begin{Ex}
{\rm We consider the quotient singularity $G=\frac{1}{5}(1,4) \subseteq SL(2,\mathbb{C})$. The group $G$ act on $\mathbb{C}^2$ by $(x,y) \mapsto (\varepsilon x,\varepsilon^4 y)$, and $\mathbb{C}^2/G$ is called $A_4$-type singularity. The lattice set $N$ is $=\mathbb{Z}^2+\sum_{\bar{g} \in \phi(G)}\mathbb{Z}\bar{g}$. Let $\rho$ be the irreducible representation of $G$ such that defined by $\rho(g)=\varepsilon^a$ for any $g=\frac{1}{5}(a,5-a) \in G$ , then ${\rm Irr}(G)$ is generated by ${\rho}$(i.e. ${\rm Irr}(G)=\{\rho^0, \rho^1, \rho^2, \rho^3, \rho^4\}$ ).\\
We write $\Gamma_i=\Gamma(I_i)$，and denote $\sigma_i$ by the cone $\sigma_i=\sigma(\Gamma_i)$. Each five boxes shows $G-{\rm graph}$ of $G$ respectively in Figure 1. The monomials outside the boxes are generaters of ideal corresponding $\Gamma_i$. 
For example, $I_1$ is an ideal generated by $(x, y^5)$，and $\Gamma_1=\{1,y,y^2,y^3,y^4\}$. In addition, a monomial $y^i$ correspond to $\rho^{5-i} \in {\rm Irr}(G)$ and $x^i$ correspond to $\rho^i$ by the map wt$:\mathcal{M} \to {\rm Irr}(G)$．Since ${\rm wt}_{\Gamma}(y^5)=1$ and  ${\rm wt}_{\Gamma}(x)=y^4$, we can obtain $\sigma_1=\{(\omega_x,\omega_y) \in N_{\mathbb{R}} \mid 5\omega_y \geq 0, \omega_x \geq 4\omega_y\}$.}

\begin{figure}[htbp]
  \begin{center}
   \includegraphics[width=100mm]{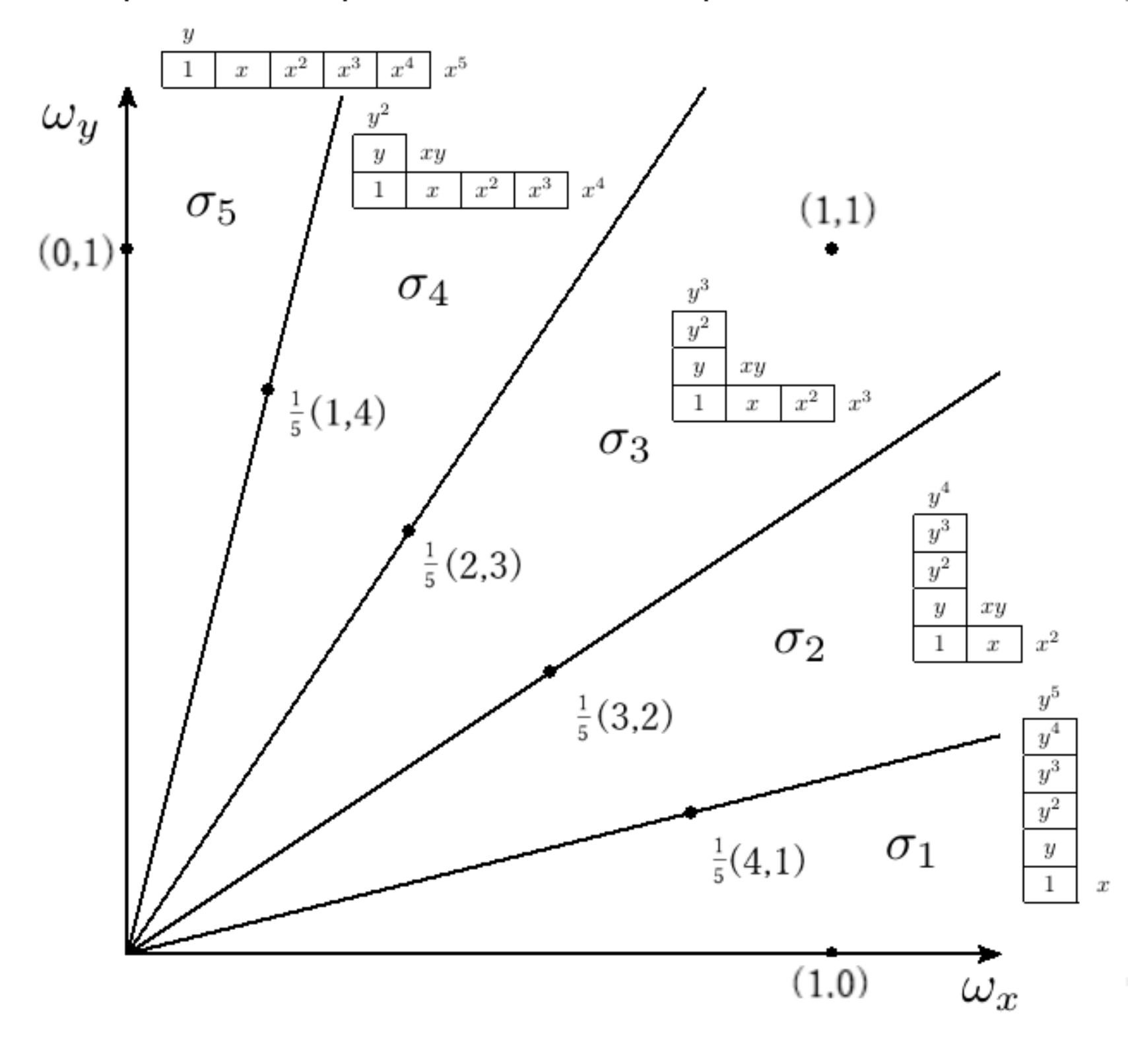}
  \end{center}
  \caption{The fan($G$) and $G${\rm -graph}}
\end{figure}
{\rm 
$I_1$ to $I_5$ are monomial ideal generated by  
$I_1=(x, y^5)$，
$I_2=(x^2, xy, y^4)$， 
$I_3=(x^3, xy, y^3)$， 
$I_4=(x^4, xy, y^2)$， 
$I_5=(x^5, y)$，respectively．Each $\sigma_i$ defined by $\Gamma_i$ is as following. }
\begin{eqnarray}
\sigma_1&=&\{(\omega_x,\omega_y) \in N_{\mathbb{R}} \mid 5\omega_y \geq 0, \omega_x \geq 4\omega_y\},  \nonumber \\
\sigma_2&=&\{(\omega_x,\omega_y) \in N_{\mathbb{R}} \mid  2\omega_x \geq 3\omega_y \geq 0, 4\omega_y \geq \omega_x\}, \nonumber \\
\sigma_3&=&\{(\omega_x,\omega_y) \in N_{\mathbb{R}} \mid  3\omega_x \geq 2\omega_y \geq 0, 3\omega_y \geq 2\omega_x\}, \nonumber \\
\sigma_4&=&\{(\omega_x,\omega_y) \in N_{\mathbb{R}} \mid  4\omega_x \geq \omega_y \geq 0, 2\omega_y \geq 3\omega_x\}, \nonumber \\
\sigma_5&=&\{(\omega_x,\omega_y) \in N_{\mathbb{R}} \mid 5\omega_x \geq 0, \omega_y \geq 4\omega_x\}. \nonumber
\end{eqnarray}

\end{Ex}

\begin{Ex}
{\rm We will show another example. We consider the cyclic quotient singularity of type $G=\frac{1}{7}(1,2,4)$ in dimension three．Lattice points of age one are $\frac{1}{7}(1,2,4)$, $\frac{1}{7}(2,4,1)$, and $\frac{1}{7}(4,1,2)$. In Figure 2, we illustrate the result of triangulation of junior simplex.
This triangulation defines a crepant resolution $Y \to X=\mathbb{C}^3/G$.} 
\\

\begin{figure}[htbp]
  \begin{center}
   \includegraphics[width=120mm]{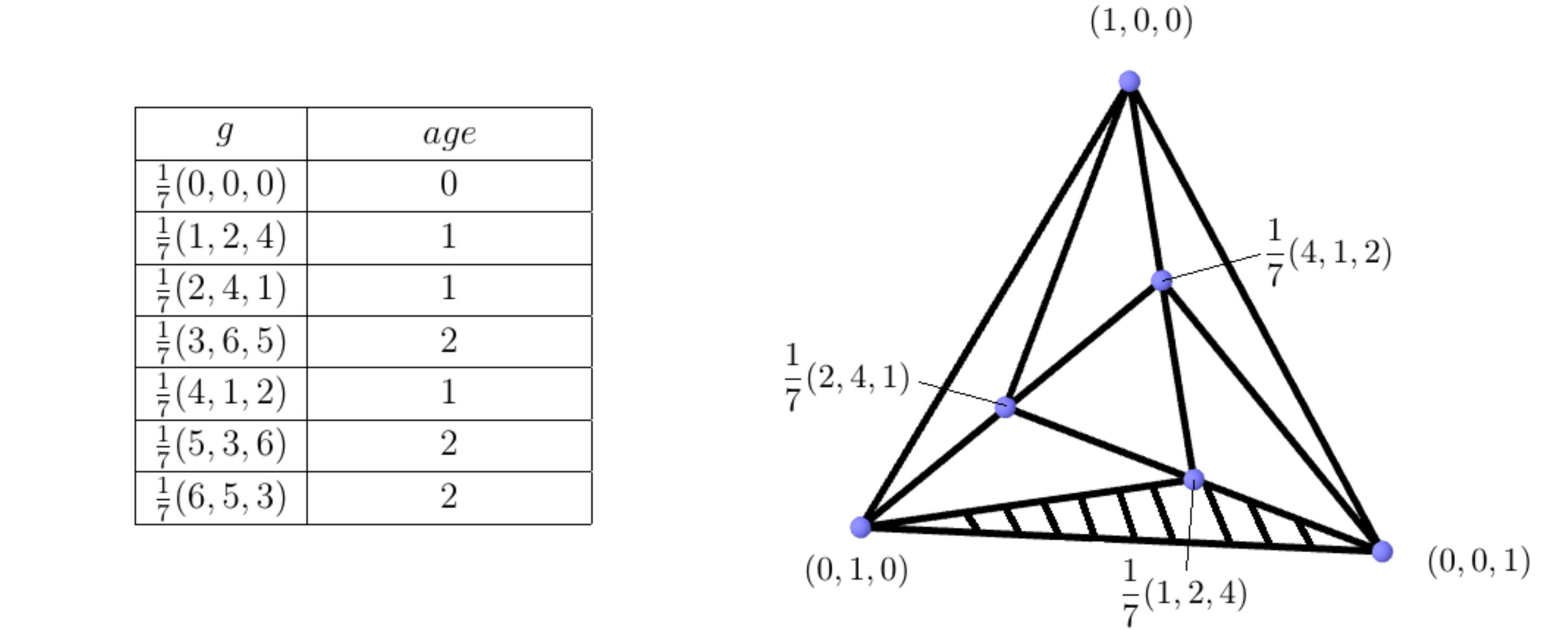}
  \end{center}
  \caption{Lattice points of age one and triangulation}
\end{figure}

\begin{figure}[htbp]
  \begin{center}
   \includegraphics[width=140mm]{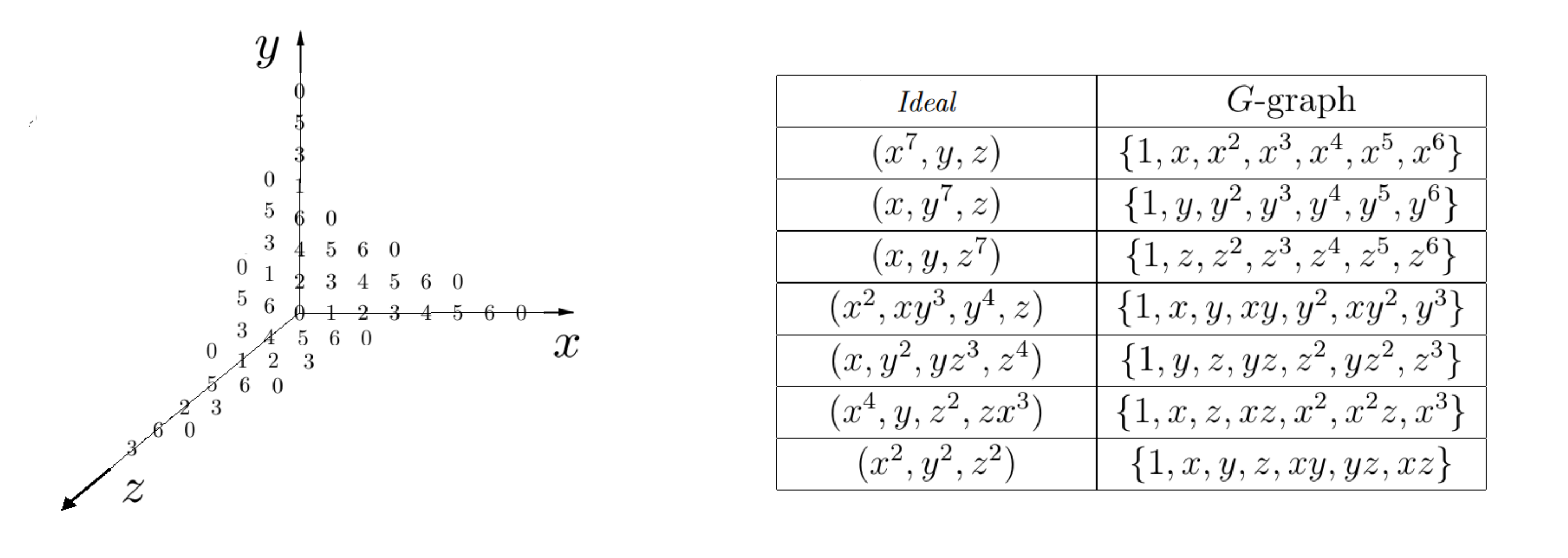}
  \end{center}
  \caption{The representation corresponding to monomials and $G$-graph}
\end{figure}

{\rm Next, we calculate  $\rm{Hilb}^G(\mathbb{C}^3)$.
As above example, the notation $\rho_i$ means that the irreducible representation of $G$ defined by $\rho_i(g)=\varepsilon^{ai}$ for any $g=\frac{1}{7}(a,b,c) \in G$, then we can find ${\rm Irr}(G)=\{\rho_0,\rho_1,\dots,\rho_6 \}$ where $\varepsilon$ is a primitive $7$th roots of unity and $\rho_0$ is trivial representation.
Figure 3 shows the relationship between monomials of $\mathbb{C}[x,y,z]$ and representations of $G$．This table means that a monomial $x^3$ correspond to the representation $\rho_3$ and $xy^2$ correspond to $\rho_5$ and so on．We can find $G$-{\rm graph} and calculate ${\rm Fan}(G)$ by this table．
For example, the $G$-{\rm graph} determined by $I=(x^7,y,z)$ is $\{1,x,x^2,x^3,x^4,x^5,x^6\}$，then $\sigma_1$ corresponding to $\Gamma_1$ is
$$
\sigma_1=\{(\omega_x,\omega_y, \omega_z) \in \mathbb{R}^3\} \mid 7\omega_x \geq 0, \omega_y \geq 2\omega_x, \omega_z \geq 4\omega_x\},
$$
since ${\rm wt}_{\Gamma}(x^7)=1$, ${\rm wt}_{\Gamma}(y)=x^2$ and ${\rm wt}_{\Gamma}(z)=x^4$ is hold.
This cone $\sigma_1$ correspond to the shaded area in Figure 2．
We calculate others $G$-{\rm graph}，then we find that the fan corresponding to $\rm{Hilb}^G(\mathbb{C}^3)$ matches the fan in Figure2．Thus, $\rm{Hilb}^G(\mathbb{C}^3)$ is crepant resolution for $\mathbb{C}^3/G$．}

\end{Ex}

\section{Result}
In this section, we will show main results of this paper. First, we introduce two series of abeliansubgroups which has crepant resolutions.
\begin{Prop}\label{main1}
If $G$ is the following type, then $\mathbb{C}^4/G$ has crepant resolutions.
\begin{itemize}
\item[{\rm(i)}] $G=\left< \frac{1}{r}(1,1,0,r-2),\frac{1}{r}(0,0,1,r-1)\right>$.
\item[{\rm(ii)}] $G=\frac{1}{r}(1,a,a^2,a^3)$, where $1+a+a^2+a^3$
\end{itemize}
\end{Prop}

For these gruops $G$, ${\rm Hilb}^G(\mathbb{C}^n)$ is not necesarily crepant resolution.
\begin{Prop}\label{main2}
If $G$ is the type (i) or (ii), the relationship between ${\rm Hilb}^G(\mathbb{C}^n)$ and $\mathbb{C}^n/G$ is as follows.
\begin{itemize}
\item[{\rm(i)}] If $G=\left< \frac{1}{r}(1,1,0,r-2),\frac{1}{r}(0,0,1,r-1)\right>$  and $r$ is odd, then ${\rm Hilb}^G(\mathbb{C}^4)$ is blow-up of certain crepant resolutions. If $r$ is odd, then ${\rm Hilb}^G(\mathbb{C}^4)$ is one of crepant resolutions. 
\item[{\rm(ii)}] If $G=\frac{1}{r}(1,a,a^2,a^3)$ and  if $a=3$, then ${\rm Hilb}^G(\mathbb{C}^4)$ is blow-up of certain crepant resolutions.
\end{itemize}
\end{Prop}

\subsection{The case of $G=\left< \frac{1}{r}(1,1,0,r-2),\frac{1}{r}(0,0,1,r-1)\right>$ }
$(1)$ $r$ is a even number\\
We assume that $r$ is a even number.
There are $\frac{r^2}{4}+r-1$ lattice points of age one which are expressed as follows.
$$\frac{1}{r}(0,0,i,r-i)\ \ where,1\leq i \leq r-1, $$
$$\frac{1}{r}(j,j,i,r-i-2j)\ \  where,1\leq j \leq \frac{r}{2} ,\ 0 \leq i \leq r-2j$$
These lattice points are on the triangle $T$ which has vertices $(0,0,0,1)$，$(0,0,1,0)$，$(\frac{2}{r},\frac{2}{r},0,0)$(Figure $4$).\\

\begin{figure}[htbp]
 \begin{minipage}{0.5\hsize}
  \begin{center}
   \includegraphics[width=70mm]{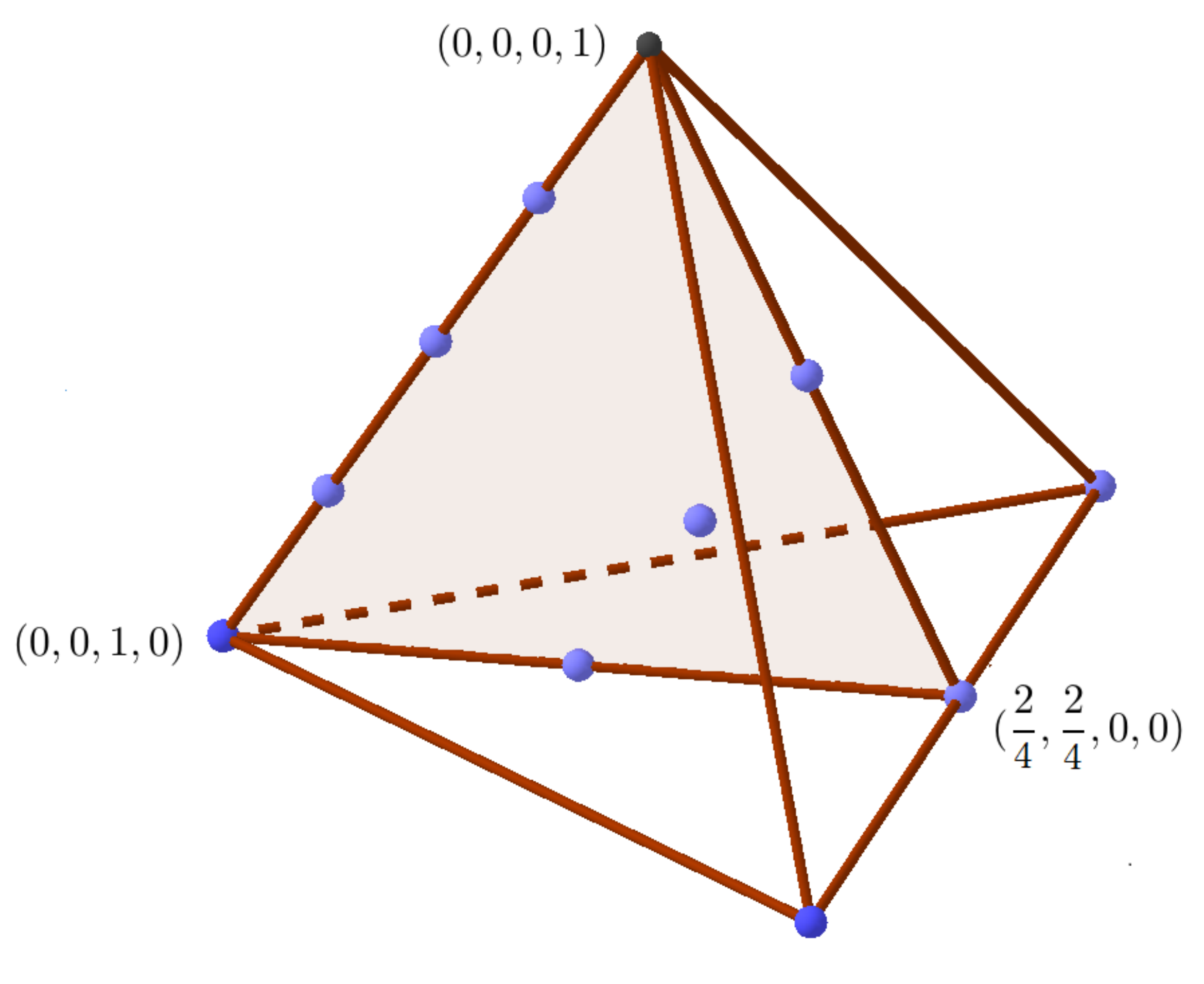}
  \end{center}
  \caption{The tetrahedron of age one}
 
 \end{minipage}
 \begin{minipage}{0.5\hsize}
  \begin{center}
   \includegraphics[width=70mm]{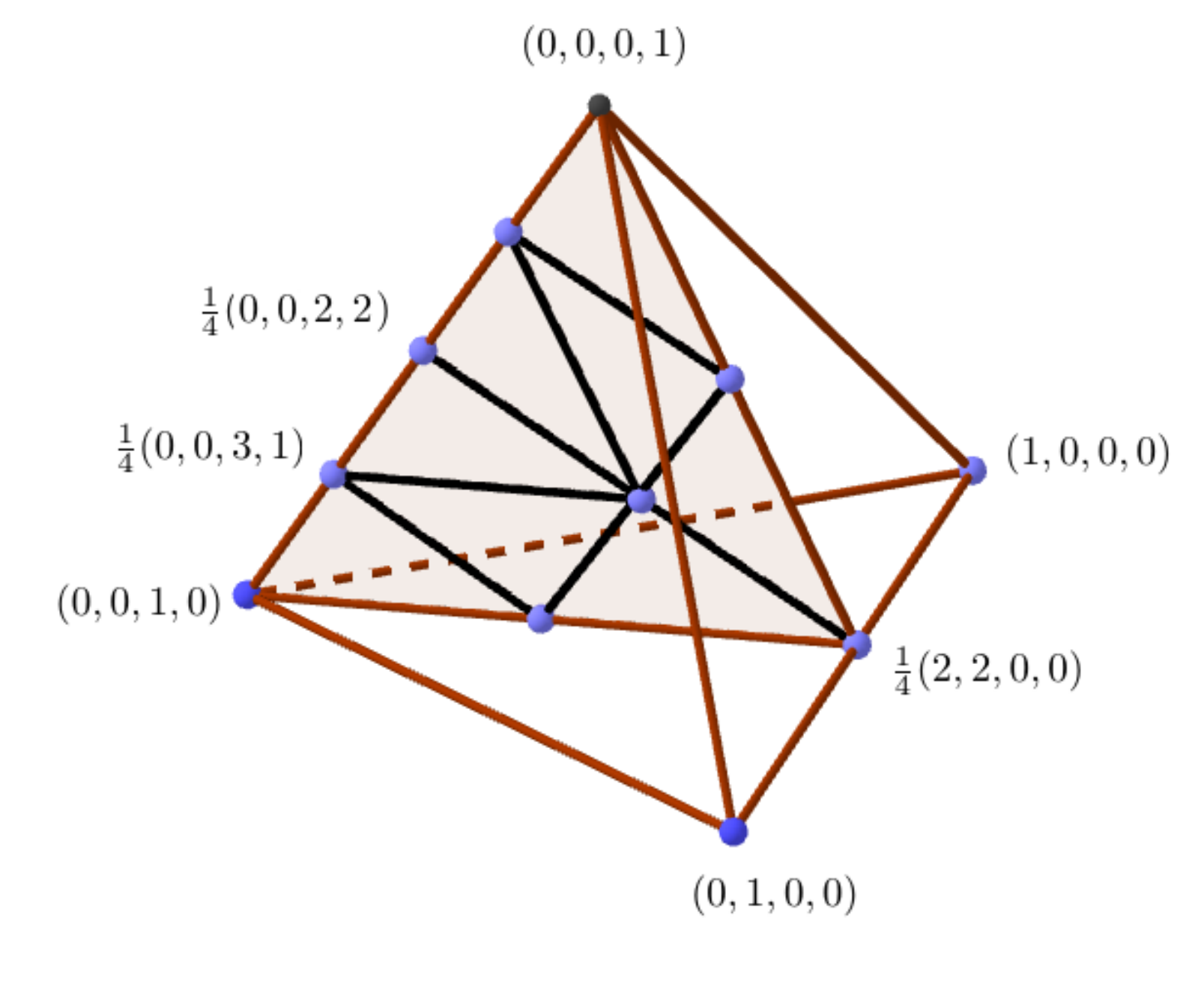}
  \end{center}
  \caption{example of triangulation}
 \end{minipage}
\end{figure}
First, we consider triangulation of $T$. After that, we choose from either $(1,0,0,0)$ or $(0,1,0,0)$ as vertex for these triangles, then we obtain tetrahedrons.
The subdivision obtained by this way gives the fan which correspond to a crepant resolution for $\mathbb{C}^4/G$. Figure $5$ show example of triangulation in the case of $r=4$. In this way, we obtain crepant resolutions as many as the number of triangulation on $T$. \\
Now for ${\rm Hilb}^G(\mathbb{C}^4)$. 
We explein about the irreducible representations of $G$ determined by monomials of  $\mathbb{C}[x,y,z,w]$ for to calculate $G$-{\rm graph}.
Arbitrary element of $G$ is  expressed  such as $\frac{1}{r}(\alpha,\alpha,\beta,\gamma)$ where $2\alpha+\beta+\gamma \equiv r \ {\rm mod}\ m$ and $m$ is a integer. Monomials $x$ and $y$ define same irreducijble reprsentations. Therefore, we write this representation $a$, and we will denote by $b$ the representation defined by $z$. We can see that monomials $x^{r-2}z^{r-1}$ and $w$ define same reprsentations $a^{r-2}b^{r-1}$. The reason is that $a^{r-2}b^{r-1}(g)=\varepsilon^{(r-2)\alpha+(r-1)\beta}=\varepsilon^{-2\alpha-\beta}$ holds for any $g=\frac{1}{r}(\alpha,\alpha,\beta,\gamma) \in G$ and a monomial $w$ correspond to a representation $a^{r-2}b^{r-1}$ since $\varepsilon^{\gamma}=\varepsilon^{rm-2\alpha-\beta}=\varepsilon^{-2\alpha-\beta}$.

In our case, there are two types of the definition ideal of $G-{\rm graph}$
\begin{eqnarray}
&type1&:\ (x^{r-2i}, y, z^{1+i+k},w^{r-k},(zw)^{1+i},z^{k+1}x^{2(k+1)},w^{r-(i+k)}x^{r-2(i+k)},x^2zw) \nonumber \\
&type2&:\ (x^{r-2i},y, z^{2+i+k},w^{r-k},(zw)^{1+i},z^{k+1}x^{2(k+1)},w^{r-(i+k+1)}x^{r-2(i+k+1)},x^2zw) \nonumber
\end{eqnarray}
where $k=0,\dots,\frac{r}{2}-1,\ i=0,\dots,\frac{r}{2}-1-k$．\\
In addition, an ideal which replaced letter $x$ and $y$ or $z$ and $w$ in above ideal defines also $G$-graph.

From now on, restrict for simplicity to the case $r=4$.
There are three ideal of type$1$ and type$2$ is only $I_2$.
\begin{eqnarray}
&I_1&=(x^4, y, z, w^4) \nonumber \\ 
&I_2&=(x^4, y, z^4, w^2, x^2w, z^3x^2, wz) \nonumber \\ 
&I_3&=(x^4, y, z^3, w^2, wz) \nonumber \\
&I_4&=(x^2, y, z^2, w^4) \nonumber
\end{eqnarray}
We need to confirm that $\Gamma(I_i)$ is $G$-graph.
The representations determined by monomials are shown in the Figure $6$. 
We can check that $\Gamma(I_i)$ is  $G$-{\rm graph} using this figure.
For example, $\Gamma(I_1)=\{1, x, x^2, x^3, w, xw, x^2w, \dots, x^2w^3, x^3w^3\}$. Since these monomials define the representation different from each other，$\Gamma(I_1)$ is $G$-{\rm graph}. Figure $7$ shows the cross section of the cone $\sigma(\Gamma(I_1))$ at age one.\\

\begin{figure}[htbp]
 \begin{minipage}{0.5\hsize}
  \begin{center}
   \includegraphics[width=72mm]{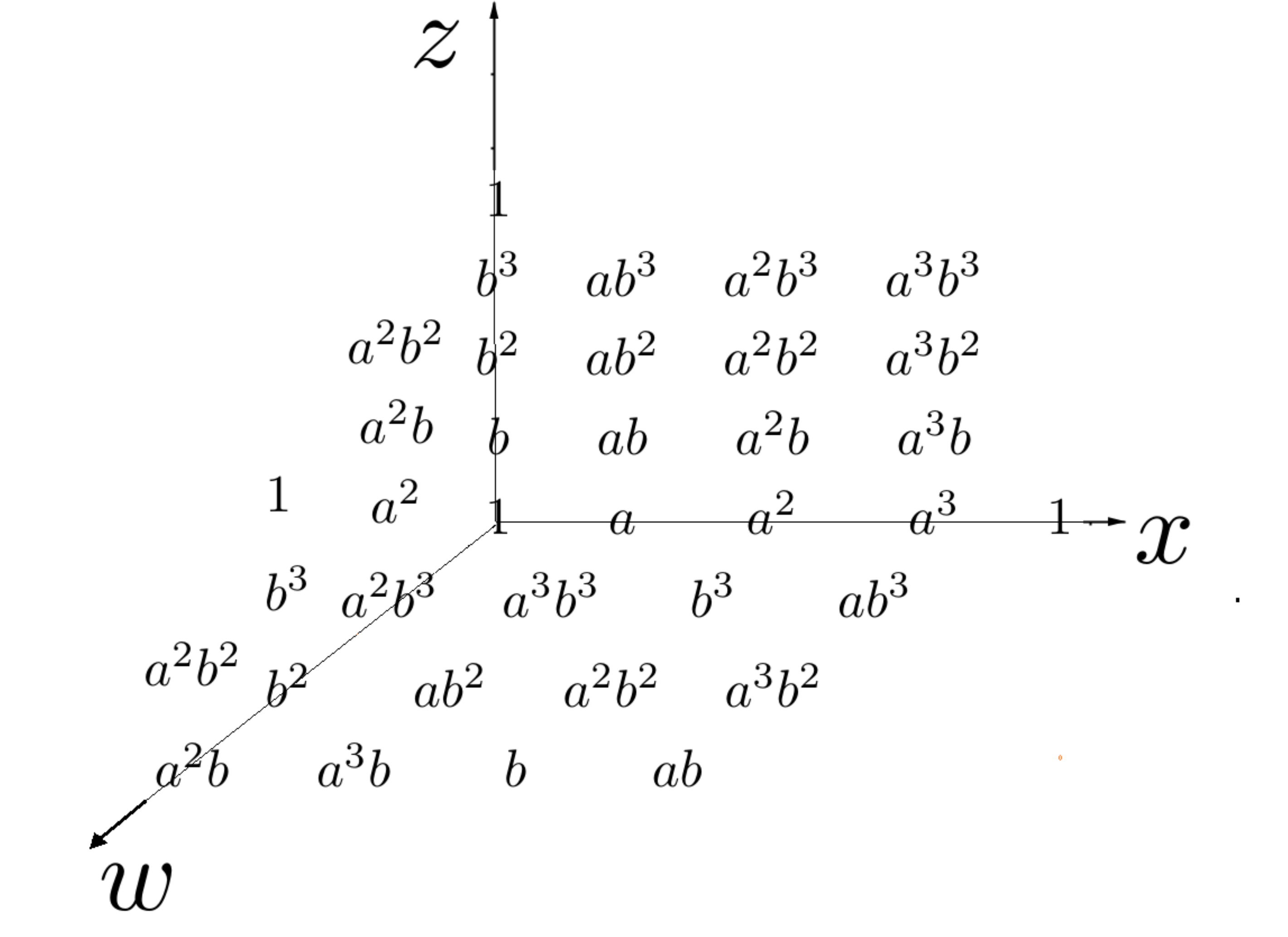}
  \end{center}
  \caption{Monomials and Representations}

 \end{minipage}
 \begin{minipage}{0.5\hsize}
  \begin{center}
   \includegraphics[width=72mm]{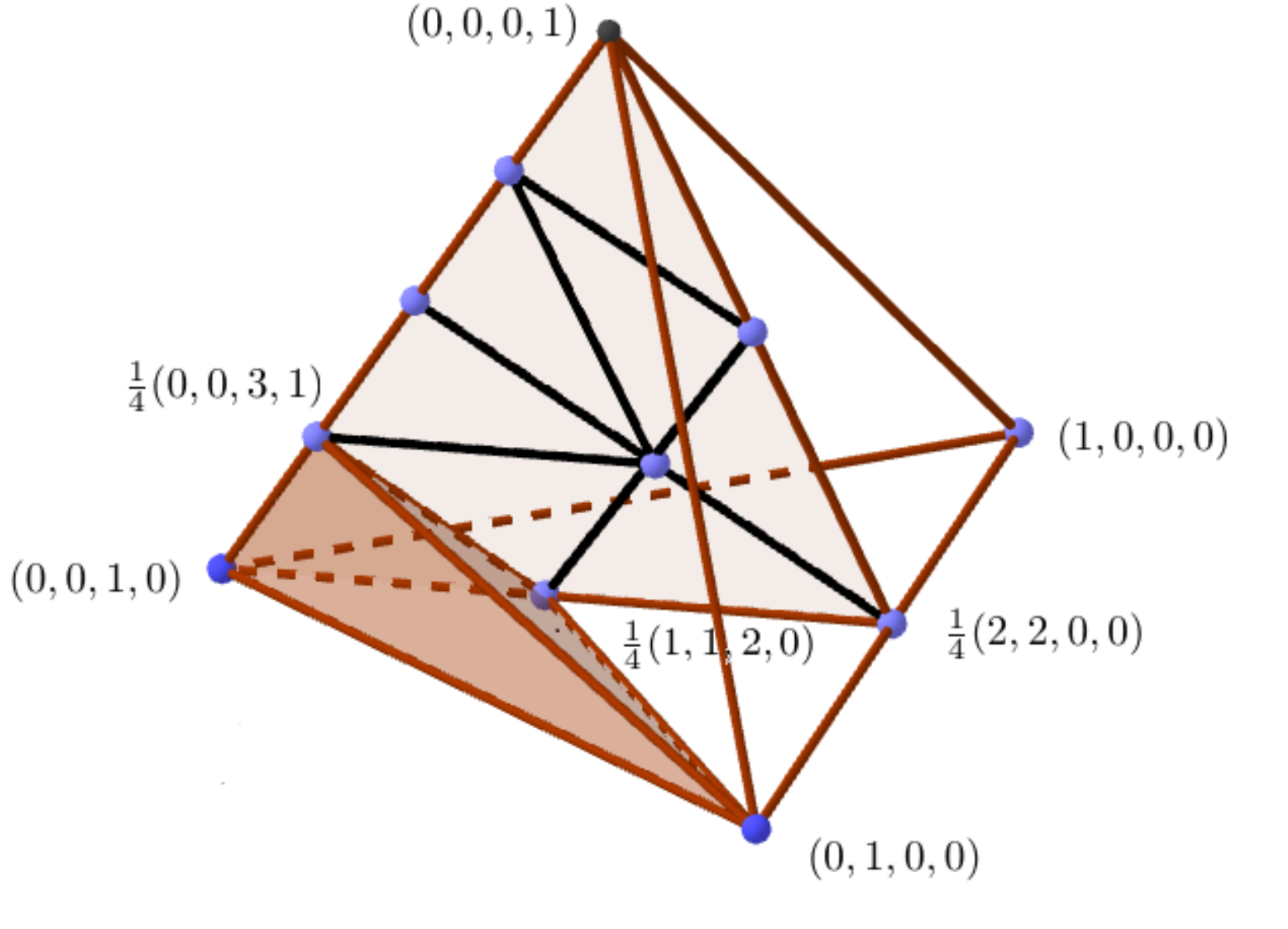}
  \end{center}
  \caption{The cone defined by $\Gamma(I_1)$}
 \end{minipage}
\end{figure}

In same way, we find that $\Gamma(I_i)$ is also $G$-graph for $I_2, I_3, I_4$
Figure $8$ shows the cone $\sigma(\Gamma(I_i))$. The difference between type$1$ and type$2$ appear in the direction of the correponding triangle(Figure $8$). 
\begin{figure}
   \begin{minipage}[b]{.31\columnwidth}
    \centering
    \includegraphics[width=\columnwidth]{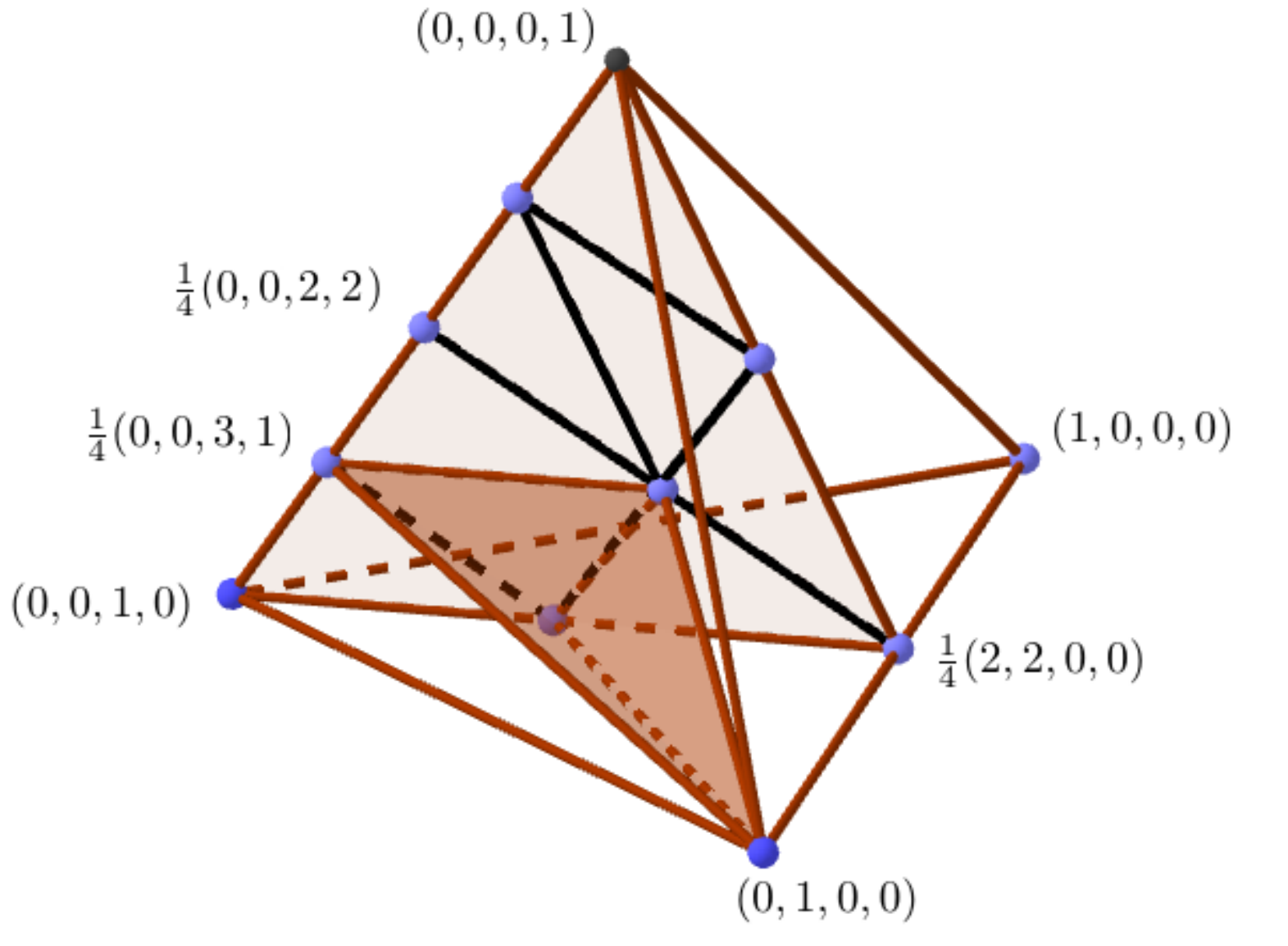} \\
    \centering $\sigma(\Gamma(I_2))$
  \end{minipage}
  \begin{minipage}[b]{.32\columnwidth}
    \centering
    \includegraphics[width=\columnwidth]{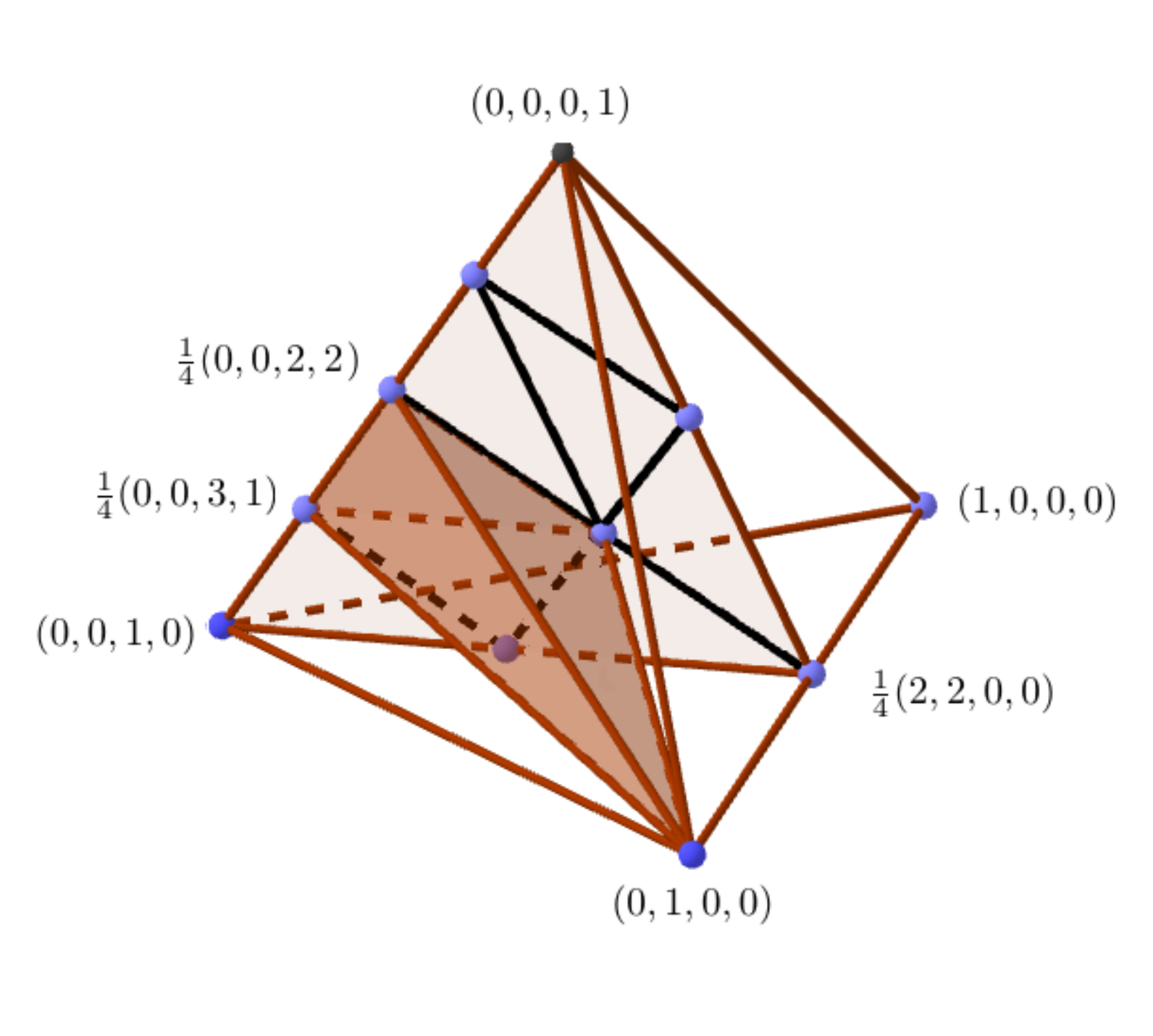} \\
       \centering $\sigma(\Gamma(I_3))$
  \end{minipage}
 \begin{minipage}[b]{.33\columnwidth}
    \centering
    \includegraphics[width=\columnwidth]{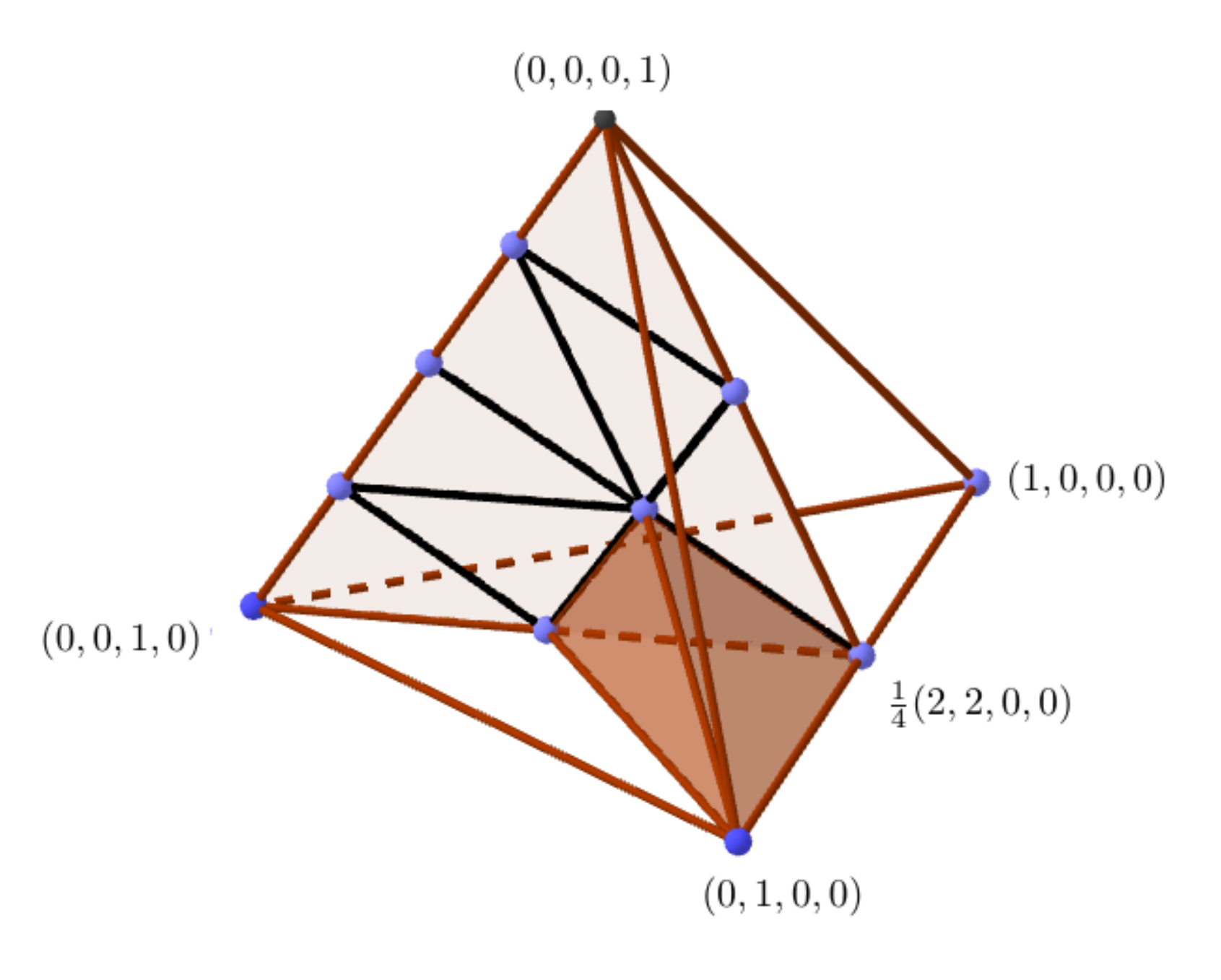}  \\
      \centering $\sigma(\Gamma(I_4))$
  \end{minipage}
 \caption{(The closs section of $\sigma(\Gamma(I_i))$)}
  \end{figure}

Others $G$-{\rm graph} are defined by ideals changing characters $z$ to $w$ and $w$ to $z$ contained in the generaters of $I_1, \dots, I_4$.We write these ideals $I_{i+4}$(where $i=1,\dots,4$).
\begin{eqnarray}
&I_5&=(x^4, y, z^4, w) \nonumber \\ 
&I_6&=(x^4, y, z^2, w^4, x^2z, x^2w^3, zw) \nonumber \\ 
&I_7&=(x^4, y, z^2, w^3, wz) \nonumber \\
&I_8&=(x^2, y, z^4, w^2) \nonumber
\end{eqnarray}
In addition, an ideal which replaced letter $x$ and $y$ in above ideal (i.e. $I_1, \dots, I_8$) defines also $G$-graph.
As a results, there are sixteen $G$-graph. Furthermore, $Fan(G)$ is shown Figure $5$

Hence, ${\rm Hilb}^G(\mathbb{C}^4)$ is a one of crepant resolutions for $\mathbb{C}^4/G$.

(2)$r$ is odd number\\
When $r$ is odd, lattice points of age one are also on the triangle $T$ , as when $r$ is even. However, the point $(\frac{2}{r},\frac{2}{r},0,0)$ is not lattice point unlike when $r$ is even (Figure 9).
\\
\begin{figure}[htbp]
 \begin{minipage}{0.5\hsize}
  \begin{center}
   \includegraphics[width=64mm]{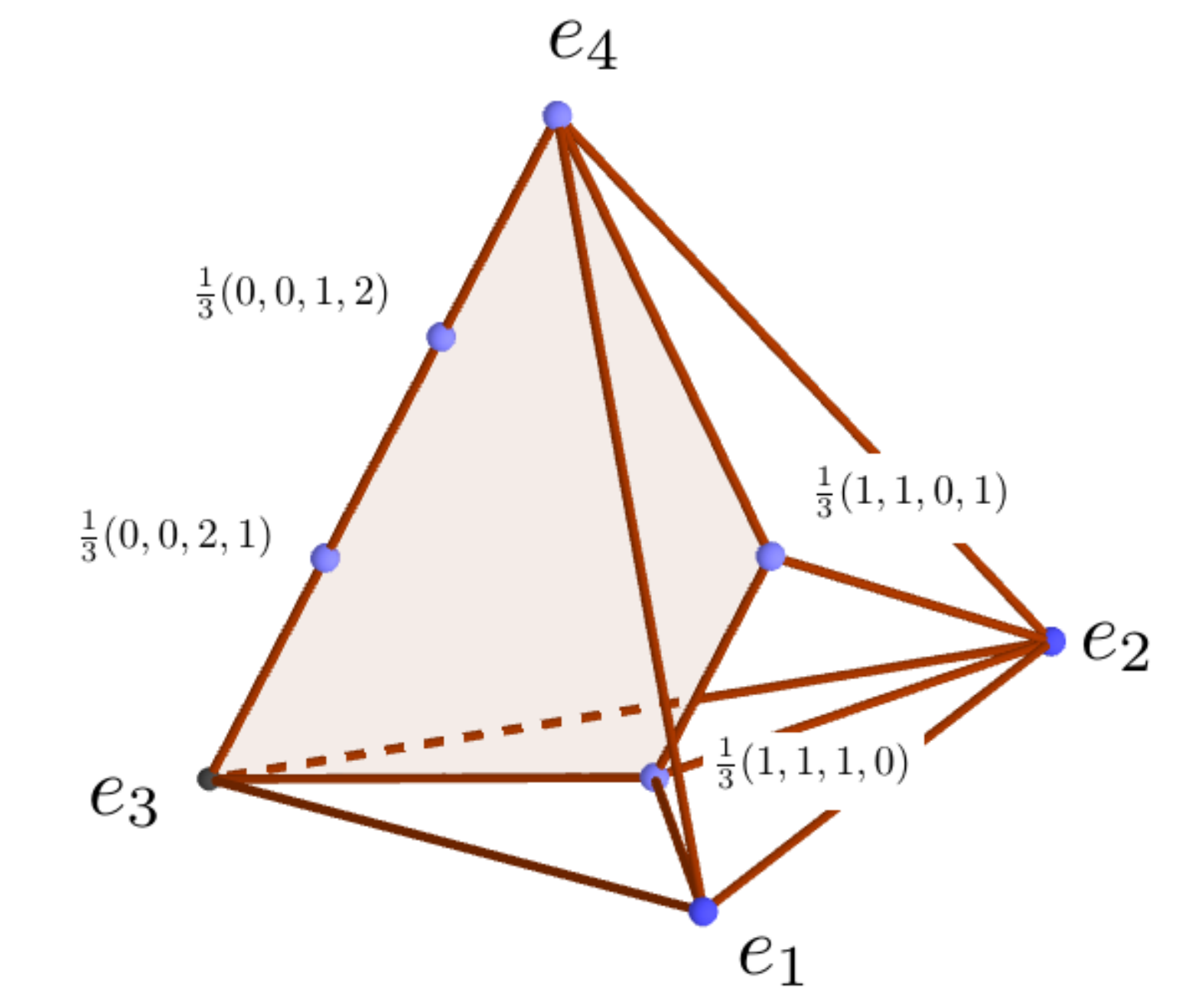}
  \end{center}
  \caption{lattice points (when $r=3$）}
 
 \end{minipage}
 \begin{minipage}{0.5\hsize}
  \begin{center}
   \includegraphics[width=58mm]{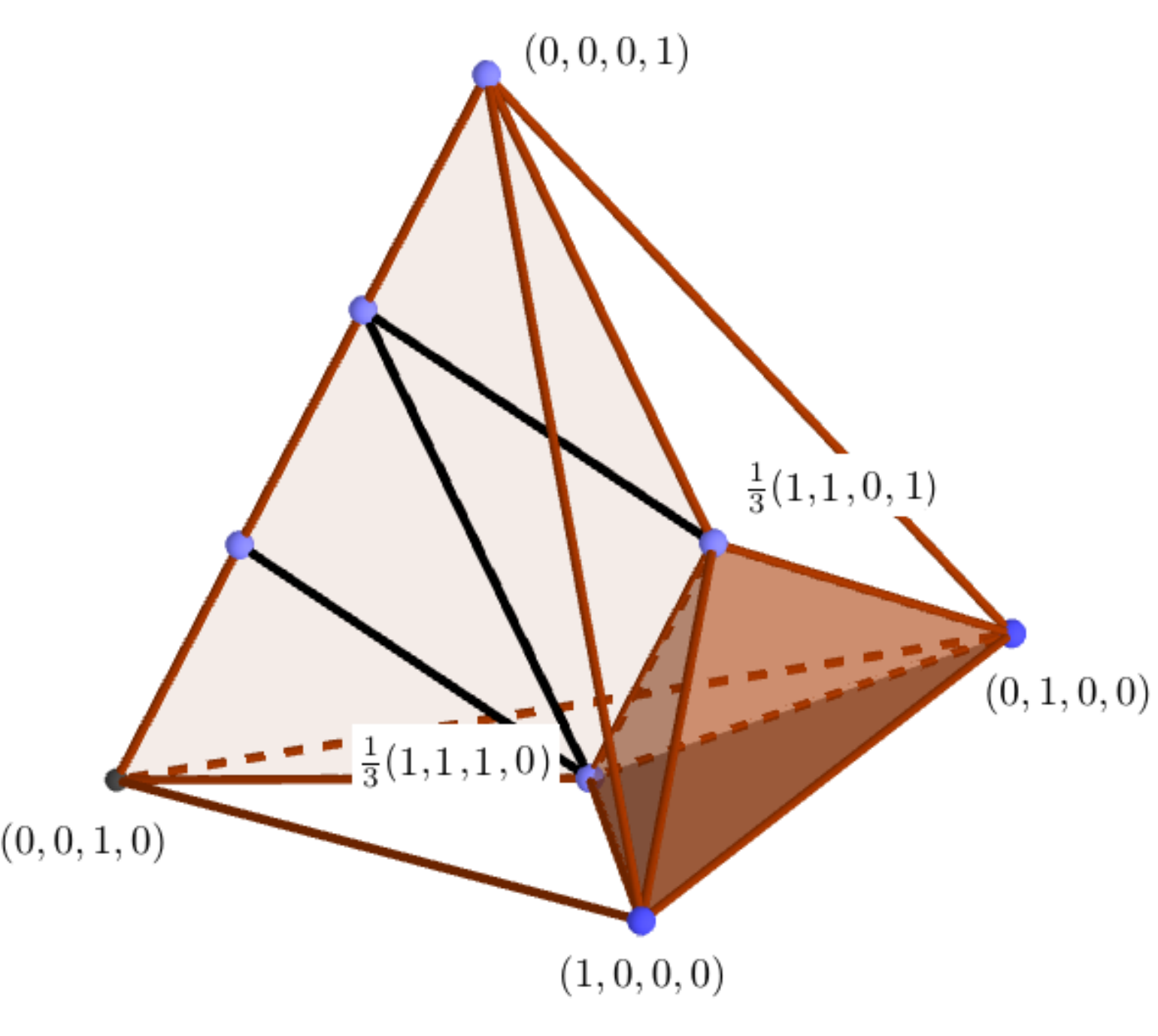}
  \end{center}
  \caption{example of crepant resolution}
 \end{minipage}
\end{figure}

Therefore, we need to consider triangulation of trapezoid instead of triangle $T$ for crepant resolution. The remain part is only a tetrahedron which has vetices $\frac{1}{r}(1,1,r-2,0)$, $\frac{1}{r}(1,1,0,r-2)$, $e_1$ and $e_2$. Figure $10$ shows example of crepant resolutions when $r=3$. \\
On the other hand, ${\rm Hilb}^G(\mathbb{C}^4)$ is not a crepant resolution.
$Fan(G)$ is the fan that subdivides the tetrahedron given by Figure $11$.
We will explain it from now.

There are six types of the definition ideal of $G-{\rm graph}$
\begin{eqnarray}
&type0:&(x,y,z^r,w^r) \nonumber \\
&type1:&(x,y^{r-2i},z^{1+i+k},w^{r-k},(zw)^{1+i},z^{k+1}y^{2(k+1)},w^{r-(i+k)}y^{r-2(i+k)},y^2zw) \nonumber \\
&type2:&(x,y^{r-2i},z^{2+i+k},w^{r-k},(zw)^{1+i},z^{k+1}y^{2(k+1)},w^{r-(i+k+1)}y^{r-2(i+k+1)},y^2zw). \nonumber
\end{eqnarray}
where $k=0,\dots,\frac{r-3}{2},\ i=0,\dots,\frac{r-3}{2}-k$．($G$-graph of type$2$ is exist when $r\geq5$）\\
\begin{eqnarray}
&type3:&\ (x,y^{r-2i},z^{\frac{r-1}{2}+i},w^{\frac{r-1}{2}},(zw)^{1+i},z^{r+1}y,w^{r+1}y,y^2zw) \nonumber \\
&type4:&\ (x,y^{r-2i},w^{\frac{r+1}{2}+i},(zw)^{1+i},z^{r+1}y,w^{r+1}y,z^{\frac{r-1}{2}-i}y^{r-1-2i},w^{\frac{r-1}{2}-i}y^{r-1-2i},y^2zw) \nonumber \\
&type5:&(x,y^{r-2i},w^{\frac{r+1}{2}+i},z^{\frac{r-1}{2}+i},(zw)^{1+i},z^{r+1}y,w^{r+1}y,z^{\frac{r-1}{2}-i}y^{r-1-2i},y^2zw) \nonumber 
\end{eqnarray}
Same as the case of even, an ideal which replaced letter $x$ and $y$ or $z$ and $w$ in above ones defines also $G$-graph. In paticular, type $0$ correspond to tetrahedron which has vetices $\frac{1}{r}(1,1,r-2,0)$, $\frac{1}{r}(1,1,0,r-2)$, $e_1$ and $e_2$.


\begin{figure}[htbp]
 \begin{minipage}{0.5\hsize}
  \begin{center}
   \includegraphics[width=60mm]{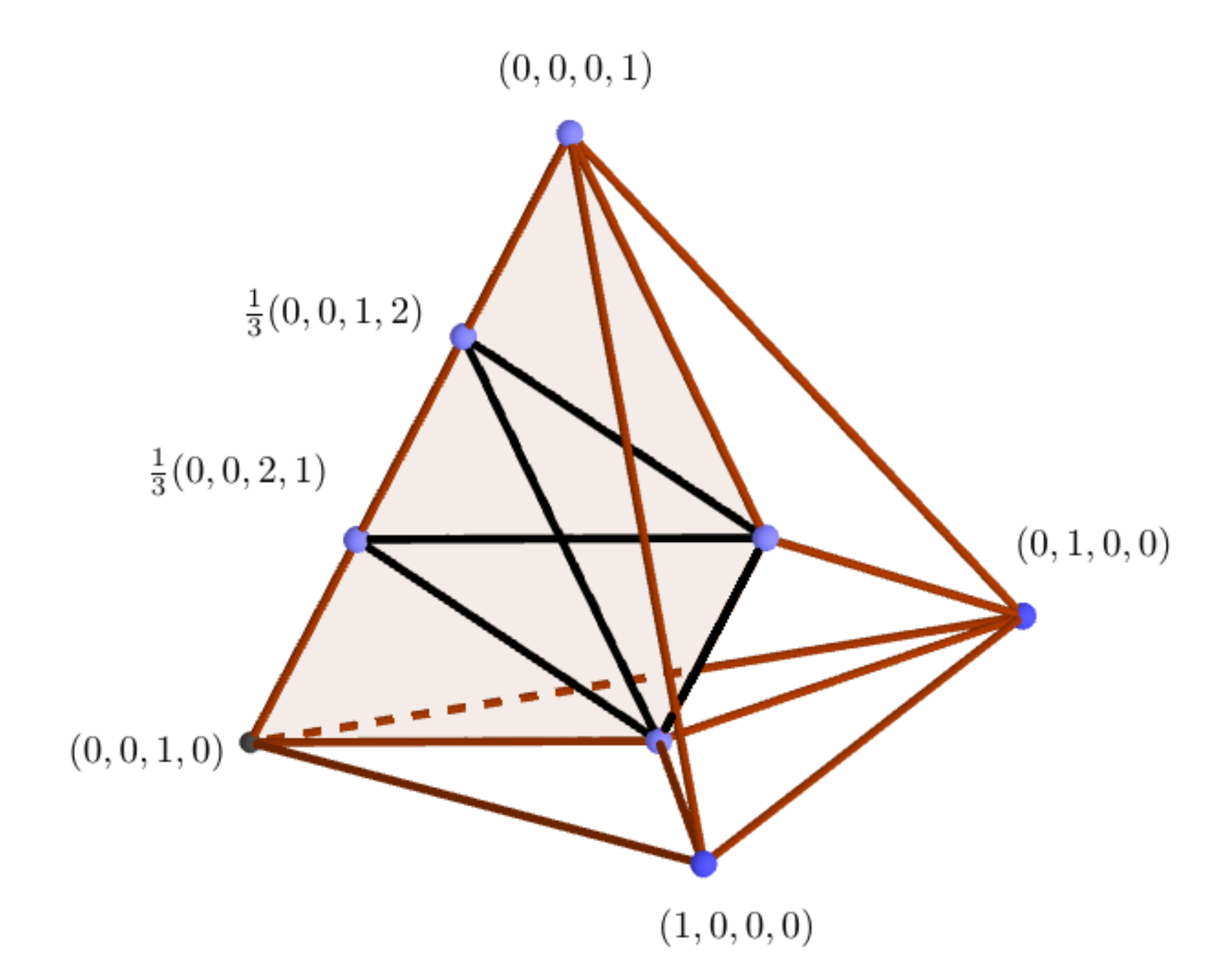}
  \end{center}
  \caption{the cross section of $Fan(G)$}
 
 \end{minipage}
 \begin{minipage}{0.5\hsize}
  \begin{center}
   \includegraphics[width=62mm]{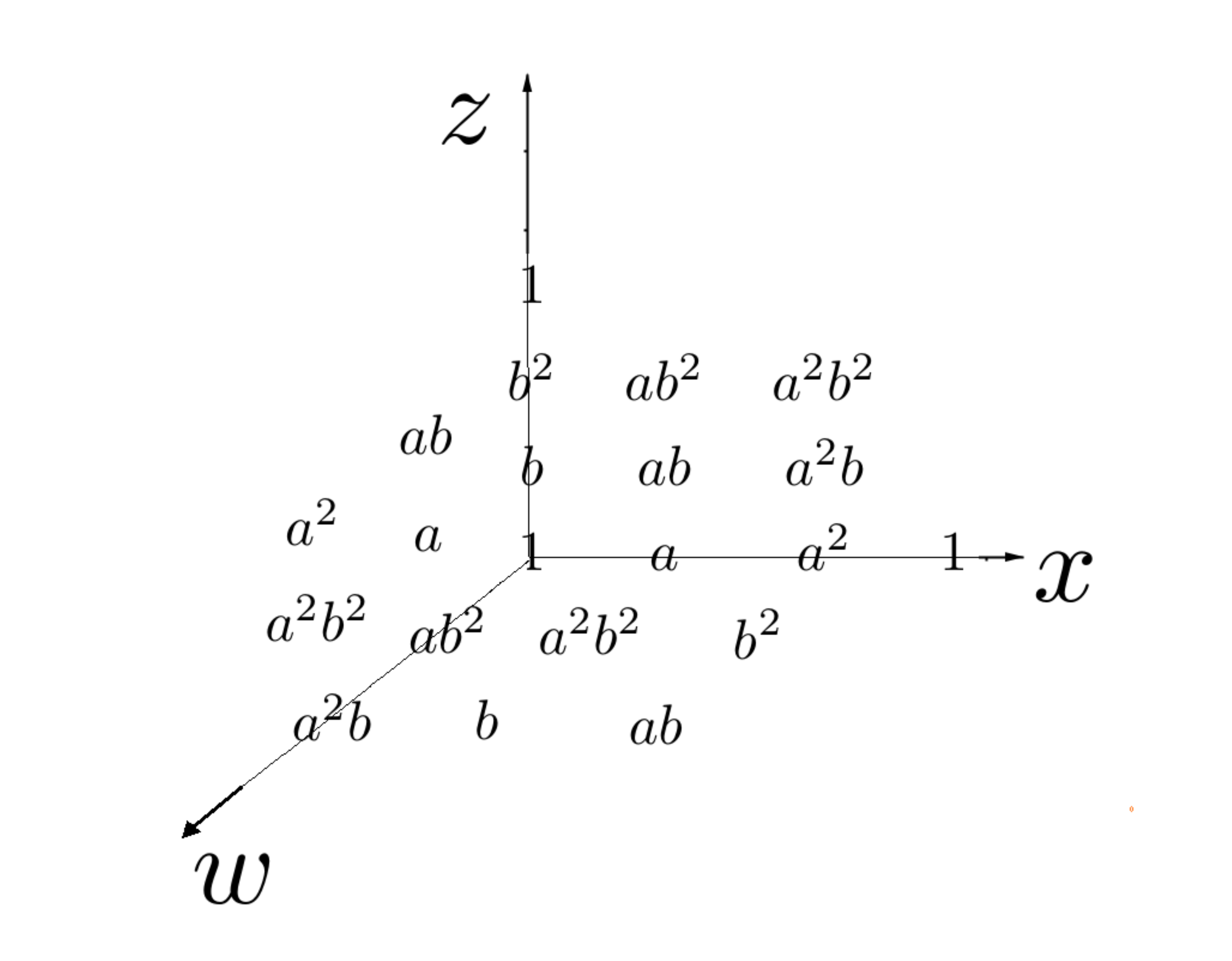}
  \end{center}
  \caption{Monomials and Representations}
 \end{minipage}
\end{figure}

For example, when $r=3$, definition ideals of $G$-{\rm graph} are
\begin{eqnarray}
type 0:&I_0&=(x, y, z^3, w^3), \nonumber \\ 
type 1:&I_1&=(x^3, y, z^3, w), \nonumber \\ 
type 3:&I_2&=(x^3 y, z^3, w^2, zw, xz^2, x^2w), \nonumber \\
type 1:&I_1'&=(x^3, y, z, w^3), \nonumber \\ 
type 3:&I_2'&=(x^3, y, z^2, w^3, zw, xw^2, x^2z), \nonumber \\
type 4:&I_3&=(x^3, y, z^3, w^3, zw, xz^2, x^2z, xw^2, x^2w), \nonumber \\
type 5:&I_4&=(x^3,y, z^2, w^2, zw). \nonumber
\end{eqnarray}
Considering replacement of letters $x$ and $y$, there are thirteen definition ideals. We can find that $\Gamma(I_i)$ is $G$-graph using Figure $12$. Figure $13$ shows the tetrahedron defined by $I_i$. 

\begin{figure}
  \begin{minipage}[b]{.24\columnwidth}
    \centering
    \includegraphics[width=\columnwidth]{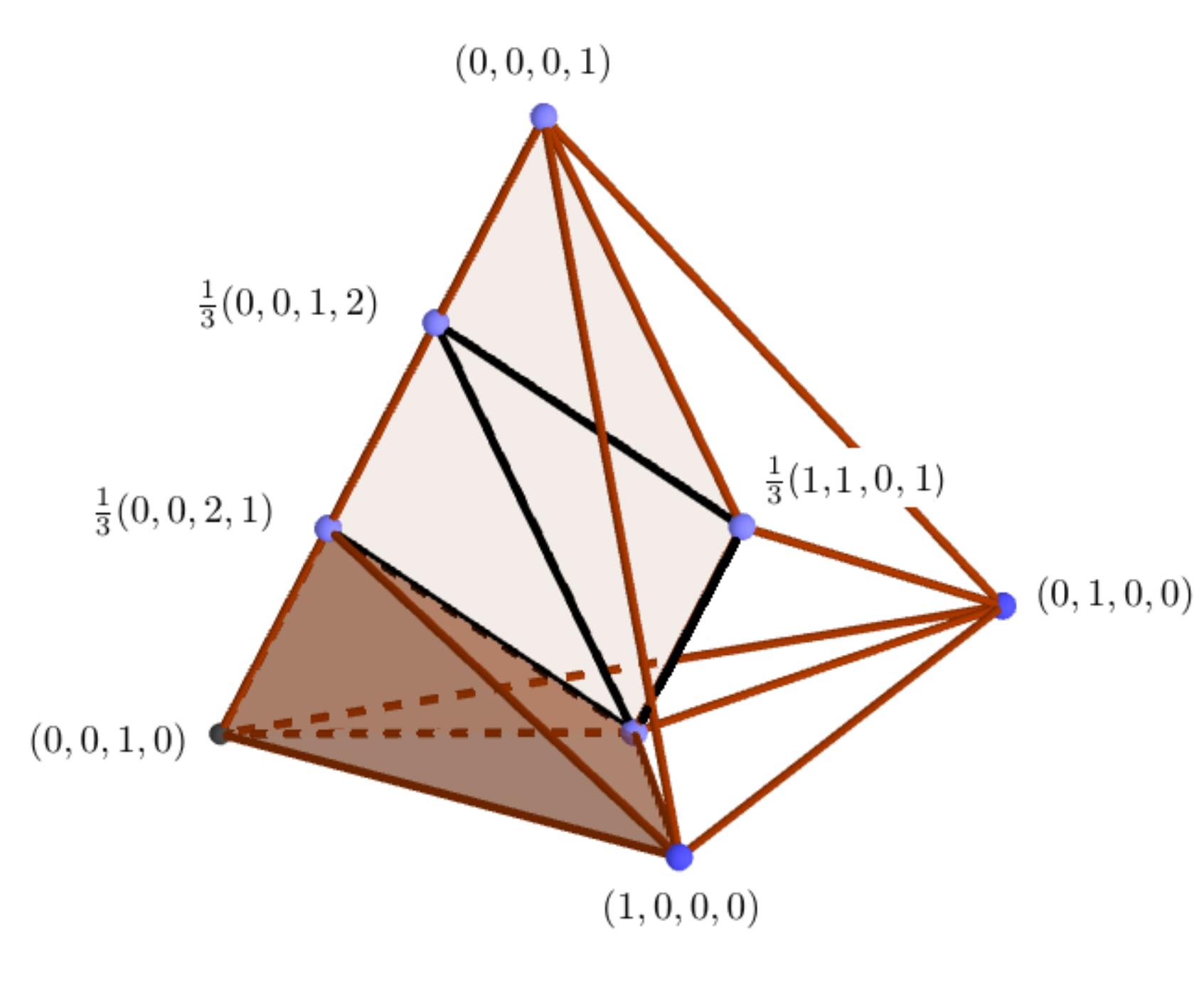} \\
    \centering $\sigma(\Gamma(I_1))$
  \end{minipage}
   \begin{minipage}[b]{.24\columnwidth}
    \centering
    \includegraphics[width=\columnwidth]{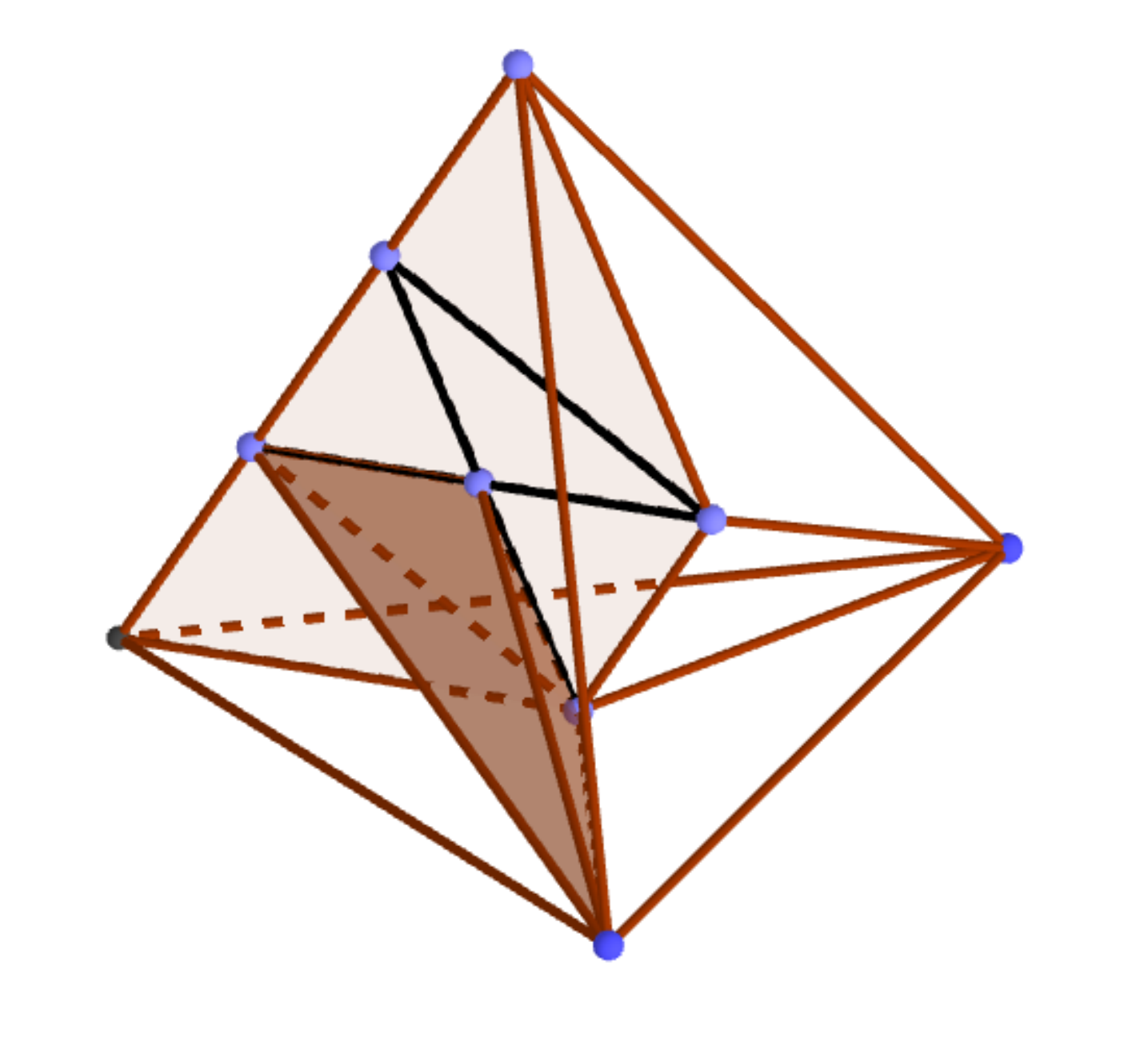} \\
    \centering $\sigma(\Gamma(I_2))$
  \end{minipage}
  \begin{minipage}[b]{.24\columnwidth}
    \centering
    \includegraphics[width=\columnwidth]{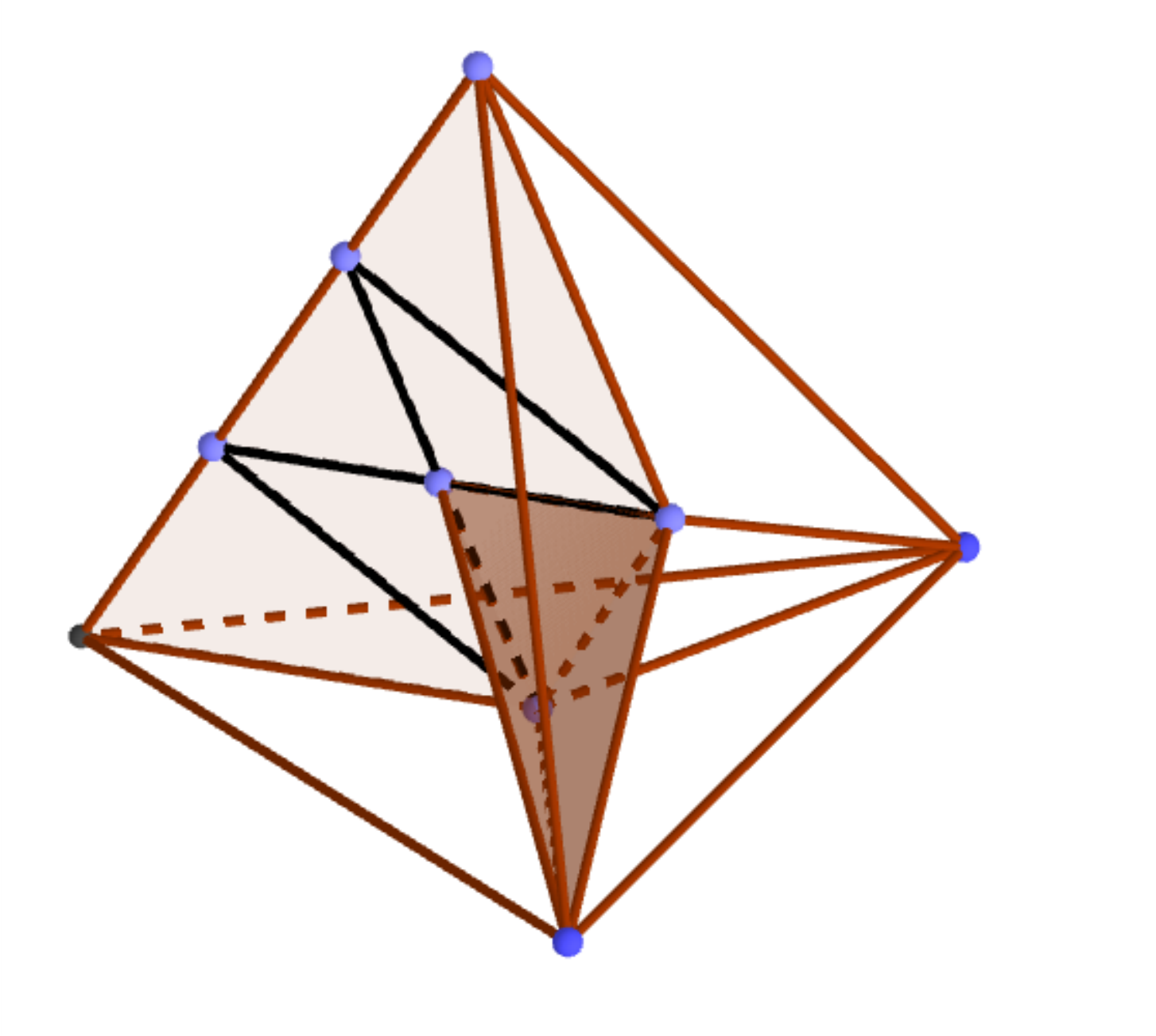} \\
       \centering $\sigma(\Gamma(I_3))$
  \end{minipage}
 \begin{minipage}[b]{.24\columnwidth}
    \centering
    \includegraphics[width=\columnwidth]{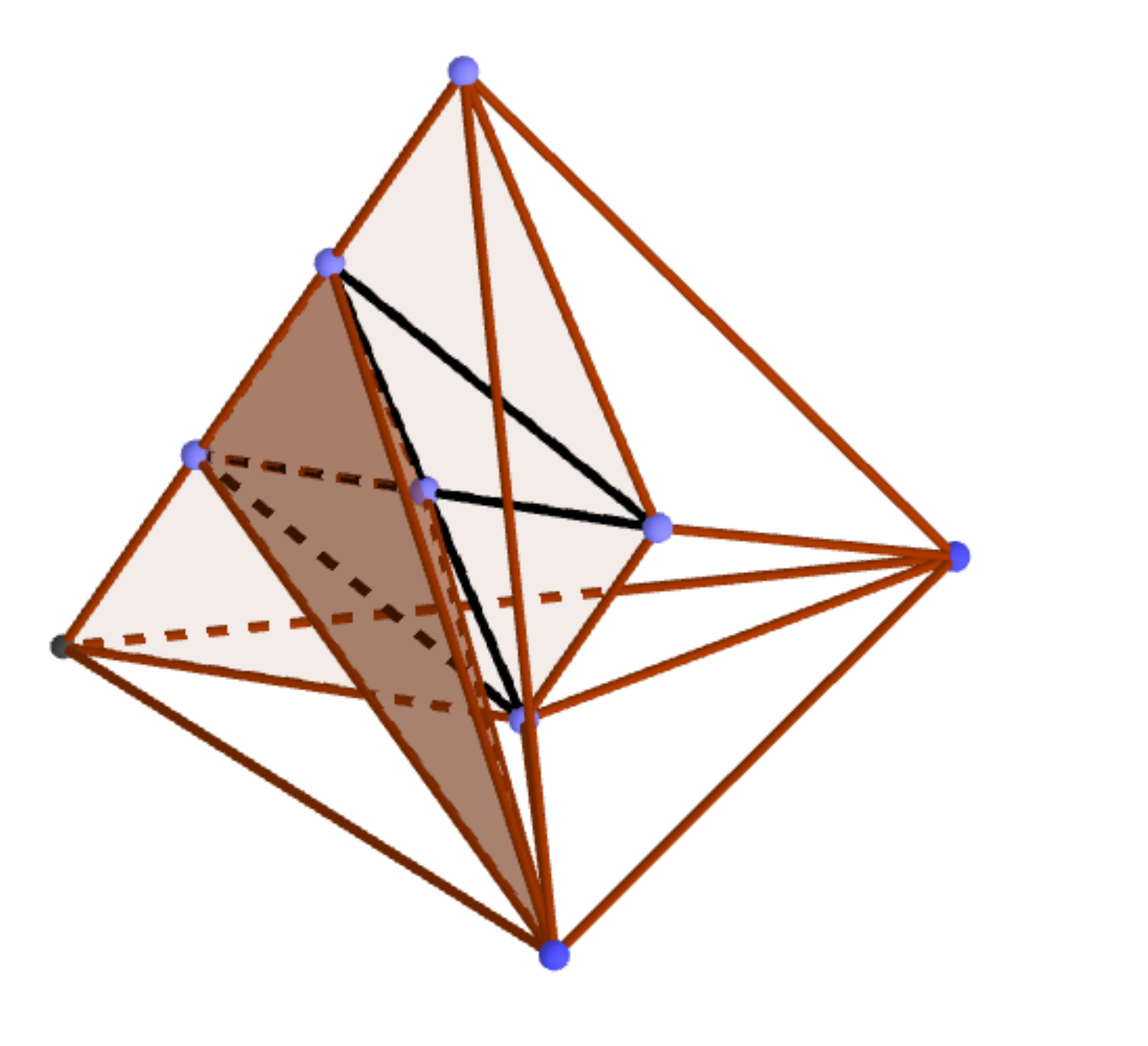}  \\
      \centering $\sigma(\Gamma(I_4))$
  \end{minipage}
 \caption{The cross section of $\sigma(\Gamma(I_i))$}
  \end{figure}

Comparing Figure $10$ and Figure $11$, we notice that $Fan(G)$ is given by subdividing the fan corresponding a crepant resolution. In other words, ${\rm Hilb}^G(\mathbb{C}^4)$ is a blow-up of one of crepant resolutions．
However, it is not a blow-up of "any" crepant resoution.
Let $Y_1$，$Y_2$ and $Y_3$ be a crepant resolution for $\mathbb{C}^4/G$(Figure $14$). Then the relationship between these and ${\rm Hilb}^G(\mathbb{C}^4)$ is as follws.
Consequently, ${\rm Hilb}^G(\mathbb{C}^4)$ is a blow-up of certain crepant resolutions.

\begin{figure}
   \begin{minipage}[b]{.31\columnwidth}
    \centering
    \includegraphics[width=\columnwidth]{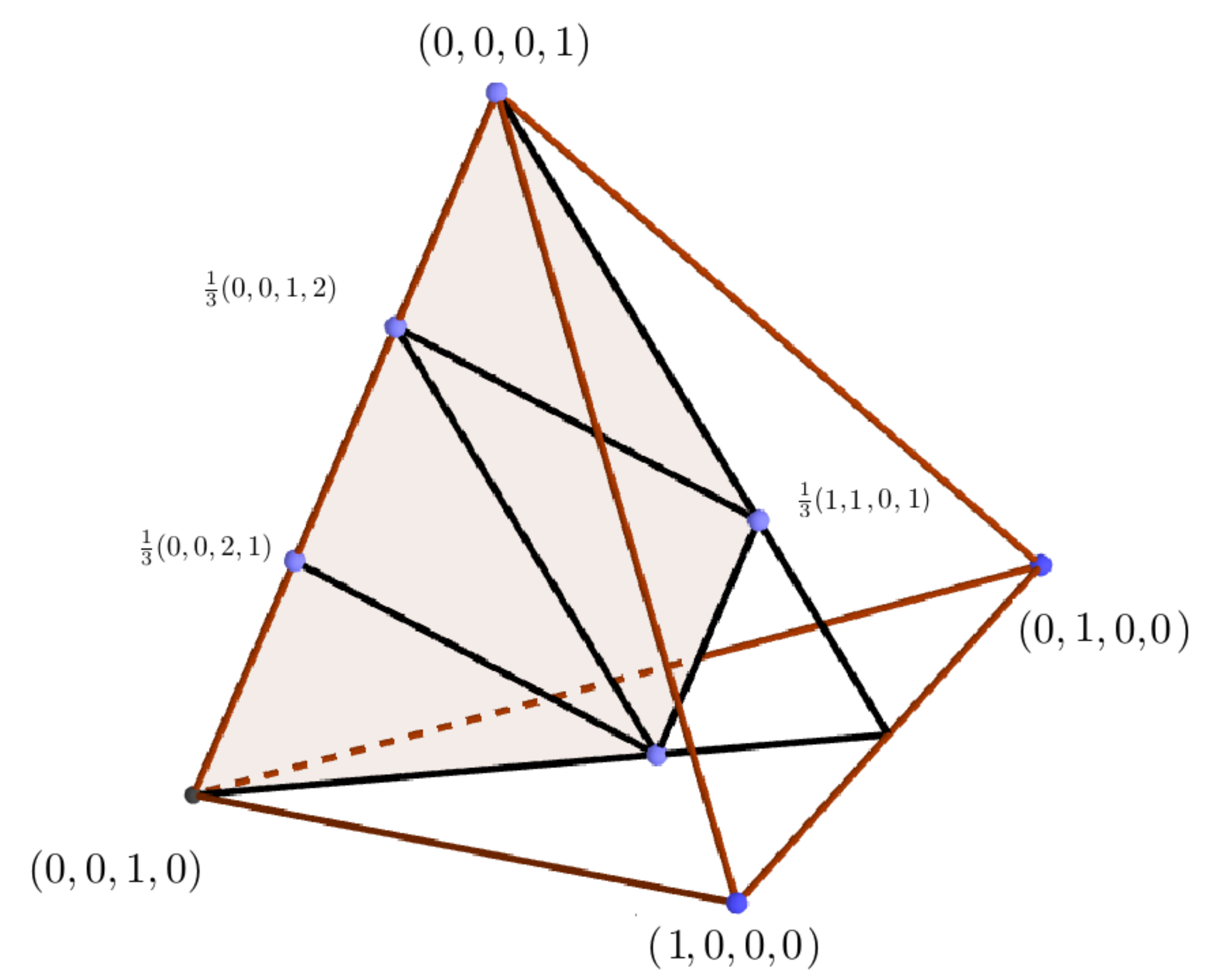} \\
    \centering $Y_1$
  \end{minipage}
  \begin{minipage}[b]{.30\columnwidth}
    \centering
    \includegraphics[width=\columnwidth]{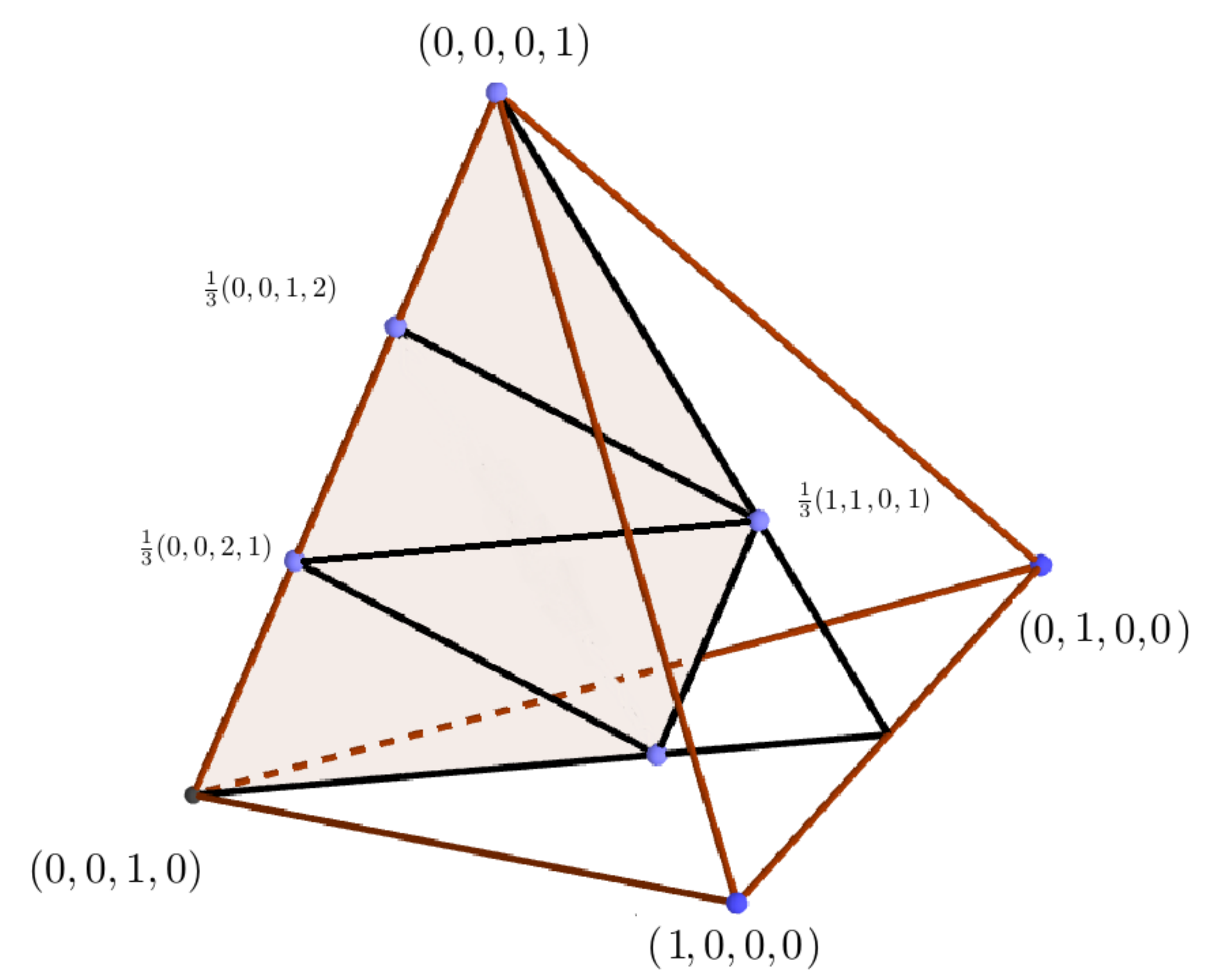} \\
       \centering $Y_2$
  \end{minipage}
 \begin{minipage}[b]{.33\columnwidth}
    \centering
    \includegraphics[width=\columnwidth]{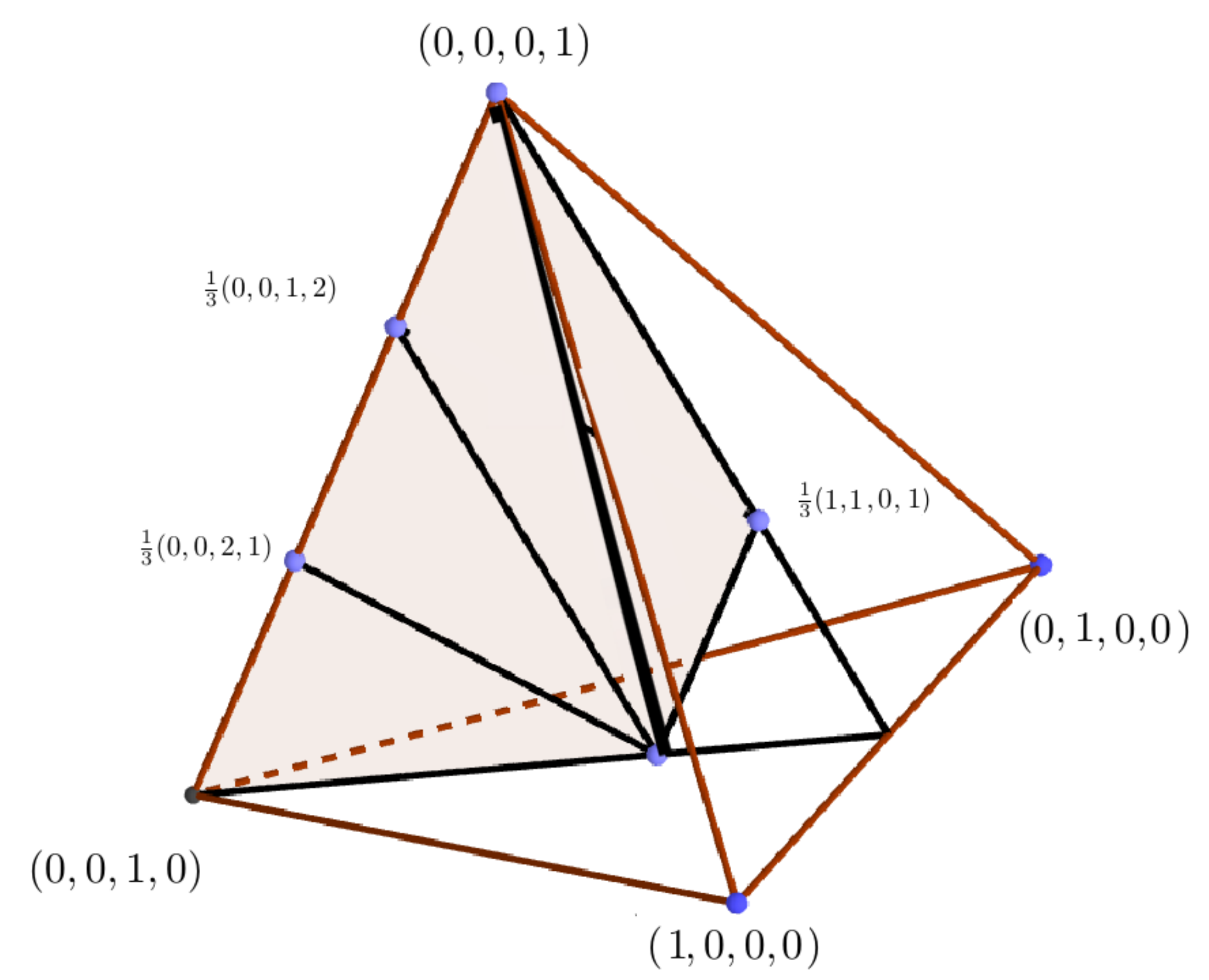}  \\
      \centering $Y_3$
  \end{minipage}
 \caption{crepant resolutions for $\mathbb{C}^n/G$}
  \end{figure}

\[
\xymatrix@R=6pt{
                    &&& Y_1 \ar[rrrd]^{f_1}              &&&  \\
 {\rm Hilb}^G(\mathbb{C}^4)\ar[rrru]|{{\rm Blow-up}} \ar[rrr]|{{\rm Blow-up}}   &&&   Y_2 \ar[rrr]^(.3){f_2}                       &&& \mathbb{C}^4/G \\
                     &&& Y_3 \ar[rrru]_{f_3}               && &
}
\]

 


\subsection{$\frac{1}{r}(1,a,a^2,a^3)$-type}
In this subsection, we will treat cyclic group $G=\frac{1}{r}(1,a,a^2,a^3)$-type where $r=1+a+a^2+a^3$. The lattice points of age one apear on the tetrahedron which has vertecis as $$
\frac{1}{r}(1,a,a^2,a^3), \frac{1}{r}(a,a^2,a^3,1), \frac{1}{r}(a^2,a^3,1,a), \frac{1}{r}(a^3,1,a,a^2)
$$ (Figure$15$). We call this tetrahedron $V$. 
\begin{figure}[htbp]
  \begin{center}
   \includegraphics[width=100mm]{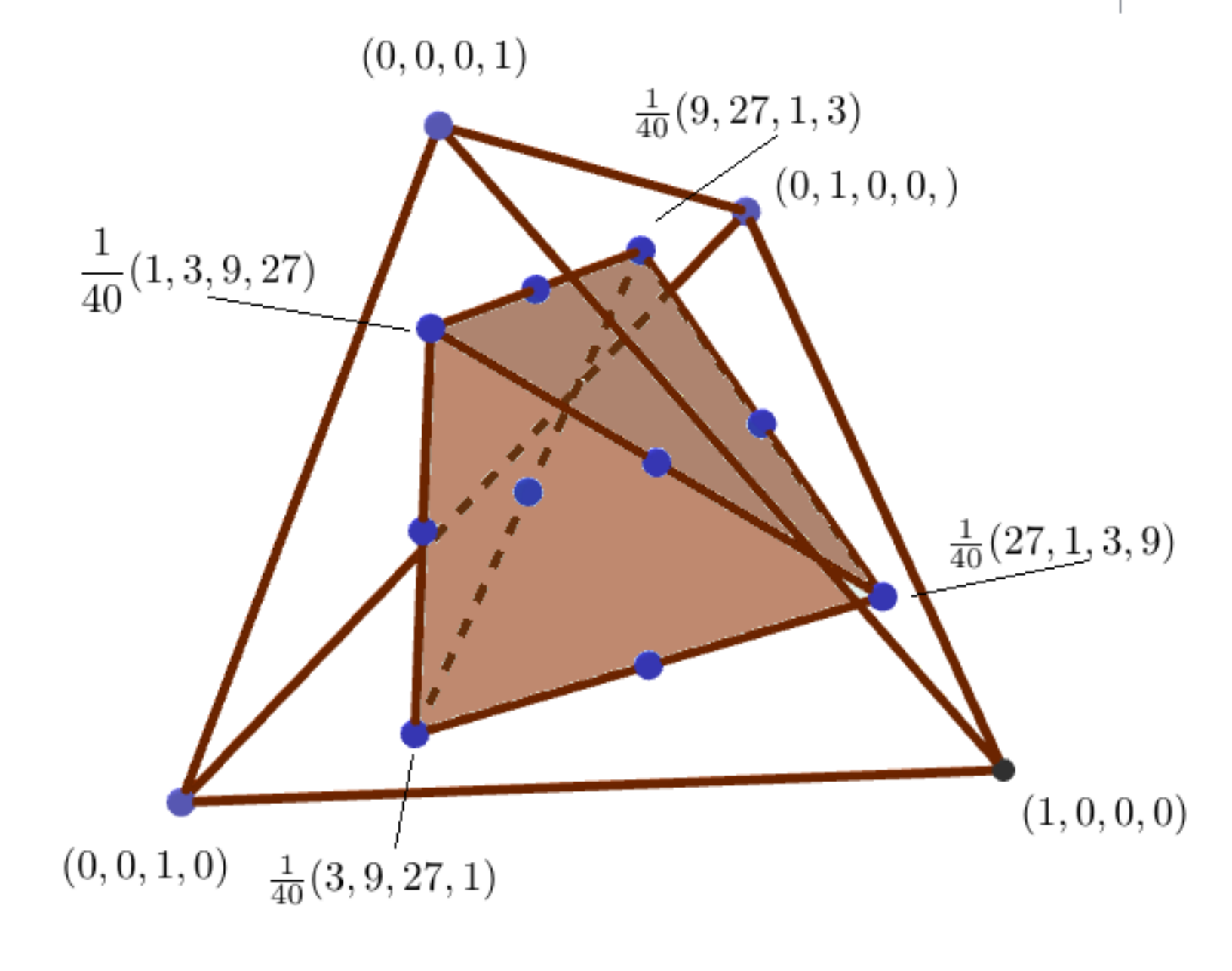}
  \end{center}
  \caption{Junior simplex when $a=3$}
\end{figure}

$V$ similer to the cross section of age one of $G'=\langle\frac{1}{a-1}(1,0,0,a-2),\frac{1}{a-1}(0,1,0,a-2),\frac{1}{a-1}(0,0,1,a-2)\rangle$.
It is well known that $\mathbb{C}^4/G'$ has crepant resolution.
We will explain tetrahedral-octahedral subdivisions of junior simplex of $\mathbb{C}^4/G'$ for crepant resolution.
First, the octahedra are obtained after chopping off tetrahedra all four couners. After that, we need to slice the remaining octahedra. We choose a pair of antipodal vertices of each octahedron and cut along this axis seem like a orange(see [11]). This method give four new tetrahedra. 
Figure $16$ and Figure $17$ shows subdivision in the case of $a=3$.
There are three way to choose the axis. Hence we can obtain three crepant resolutions in this way.
If $a>3$, then more octahedral appears in junior simplex. 

\begin{figure}[htbp]
 \begin{minipage}{0.5\hsize}
  \begin{center}
   \includegraphics[width=65mm]{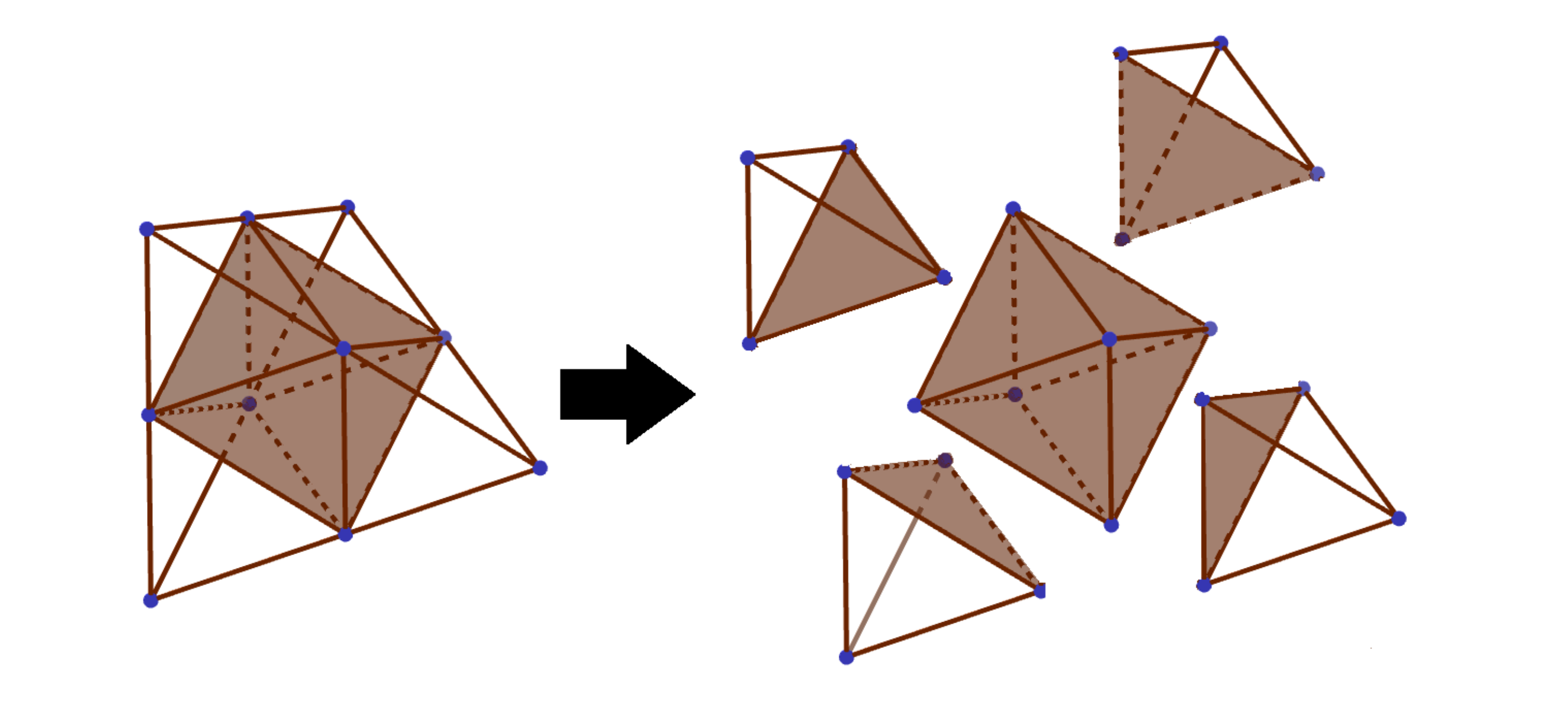}
  \end{center}
  \caption{the division of $V$}
 
 \end{minipage}
 \begin{minipage}{0.5\hsize}
  \begin{center}
   \includegraphics[width=60mm]{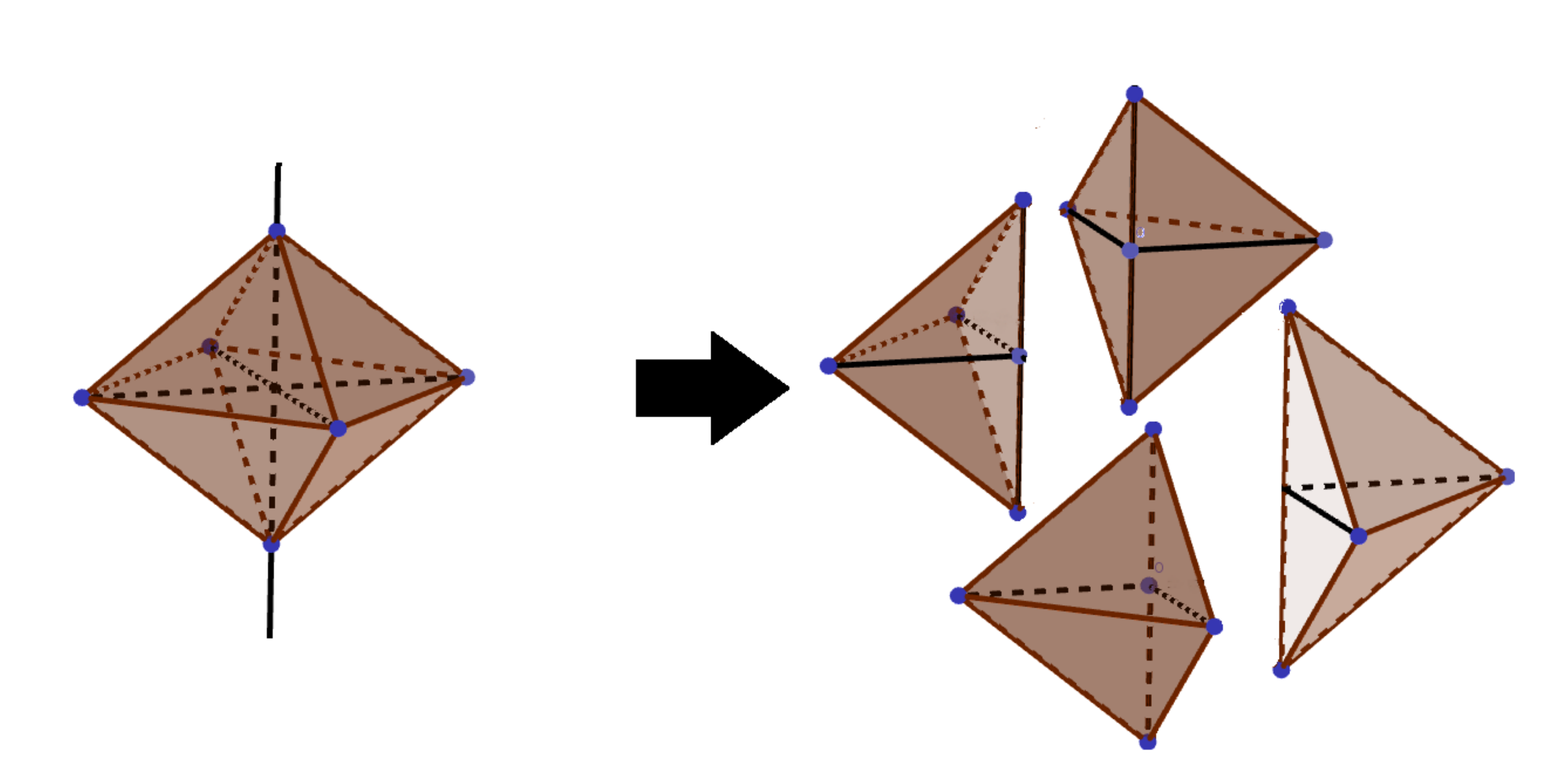}
  \end{center}
  \caption{orange-slice}
 \end{minipage}
\end{figure}
In the above argument, $V$ is divided into $(a-1)^3$ cone.
Next, we consider subdivision the outside of $V$. We call the pyramid with $e_1, e_2, e_3$ and $e_4$ at the vertex as $V_{\sigma}$.
There are 3-types tetrahedron. These tetrahedra share face, vertex or edge with $V_{\sigma}$ respectively, and we call these tetrahedra face-pyramids, vertex-pyramids and edge-pyramids.
There are four face-pyramids determined by each face $e_1, e_2, e_3$ and $e_4$. 
In addition, in the case of edge-pyramids, $a-1$ tetrahedra are obtained for each edge, so there are $6(a-1)$ tetrahedra. 
Also, there are $(a-1)^2$ vertex-pyramids for each vertex. Thus the number of vertex-pyramids is $4(a-1)^2$.

Figure $18$ shows the tetrahedron of each type.
The left picture shows face-pyramid, the middle two shows vertex-pyramid, the remaining show edge-pyramid in the case of $r=3$.

\begin{figure}[htbp]
  \begin{center}
   \includegraphics[width=140mm]{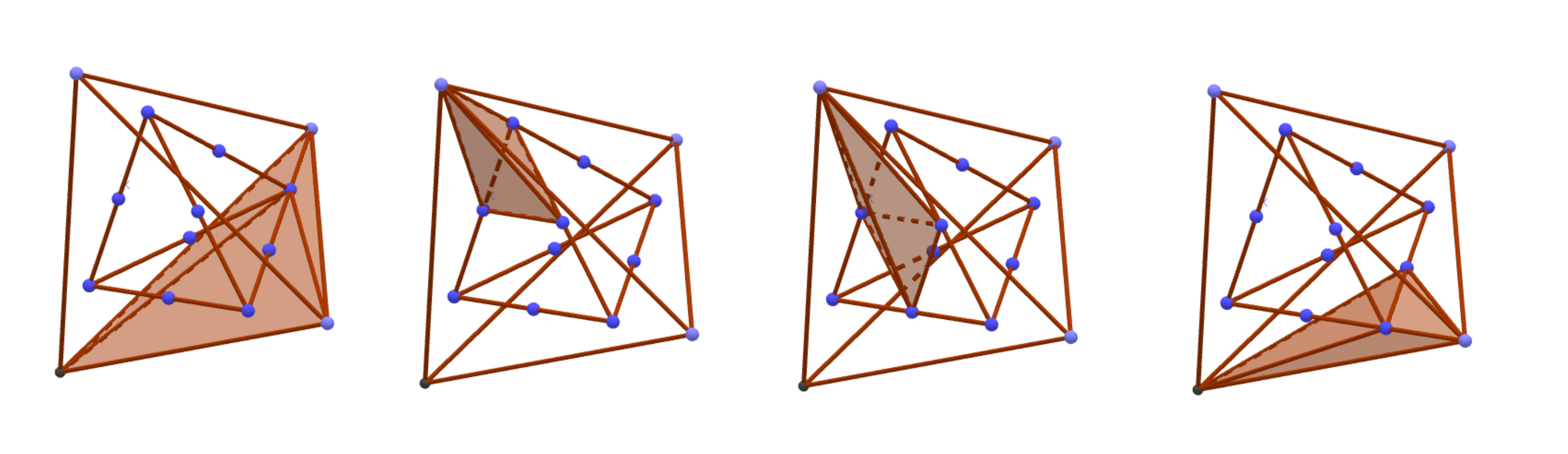}
  \end{center}
  \caption{Subdivision of outside of $V$}
\end{figure}



We recall that the number of tetrahedra in $V$ equals $(a-1)^3$.
Therefore, the total number of tetrahedra is $4+6(a-1)+4(a-1)^2+(a-1)^3=1+a+a^2+a^3=r$ in this subdivision.
The smoothness of each cone can be checked easily, so we obtain crepant resolution.



On the other hand, in the case of $a=3$, ${\rm Hilb}^G(\mathbb{C}^4)$ is not a crepant resolution, but a blow-up of it.
The following ideals define $G$-graph correspondeing to outside tetrahedron of $V$.
\begin{eqnarray}
 face-pyramids&:& (x^{40}, y, z, w) \nonumber \\
 vertex-pyramids&:& (x^3, y^3, z^5, w, z^4y^2, z^4xy) \nonumber \\
       &:& (x^3, y^3, z^9, w, z^4xy, z^5x) \nonumber \\
       &:& (x^3, y^3, z^9, w, x^2y,  z^4y^2, z^4xy, z^5x) \nonumber \\
       &:& (x^3, y^3, z^9, w, x^2y, z^4xy) \nonumber \\
 edge-pyramids&:& (x^3, y^{14}, z, w, xy^{13}) \nonumber \\
       &:& (x^2,y^{27}, z, w, xy^{13}) \nonumber
\end{eqnarray}
In paticular, Figure $18$ shows tetrahedron defined by
$$
(x^{40}, y, z, w),(x^3, y^3, z^5, w, z^4y^2, z^4xy),
(x^3, y^3, z^9, w, z^4xy, z^5x),
 (x^3, y^{14}, z, w, xy^{13}),
$$
from left to raight.\\
By replacing the letter of generater, we obtain new definition ideals of $G$-graph(for example $(x,y^4,z,w))$).
The remaining $G$-graph is given by the follwing ideals.
\begin{eqnarray}
&type\alpha&: (x^3, y^3, z^3, w^2, z^2w, y^2zw, z^2xw, xyzw) \nonumber \\
&type\beta&: (x^3, y^3, z^3, w^3, z^2w, y^2zw, z^2xw, xyzw, x^2yz, w^2y^2, w^2z^2, w^2x)\nonumber \\
&type\gamma&: (x^3, y^3, z^3, w^3, z^2w, y^2zw, z^2xw, xyzw, x^2yz, x^2z^2, w^2z^2, w^2x)\nonumber 
\end{eqnarray}
These $G$-graph define tetrahedron in $V$.
Type $\alpha$ correspond to four couners tetrahedron in tetrahedra-octahedra subdivision of $V$. We consider four tetrahedra given by subdivision of octahedron.
We cut these tetrahedra into exact halves at the midpoint of antipodal verticies, then we obtain eight new tetrahedra. Type $\beta$ and type $\gamma$ correspond to these eight tetrahedra. Theregore, ${\rm Hilb}^G(\mathbb{C}^4)$ consist of $44$ tetrahedra. Since $\mid G \mid$, ${\rm Hilb}^G(\mathbb{C}^4)$ is not a crepant resolution.

The relationship between ${\rm Hilb}^G(\mathbb{C}^4)$ and crepant resolutions $Y_1, Y_2, Y_3$ which are obtained by the tetrahedra-octahedra subdivision is as follows.

\[
\xymatrix@R=6pt{
                    &&& Y_1 \ar[rrrd]^{f_1}              &&&  \\
 {\rm Hilb}^G(\mathbb{C}^4)\ar[rrru]|{{\rm Blow-up}} \ar[rrrd]|{\rm {Blow-up}} \ar[rrr]|{{\rm Blow-up}}   &&&   Y_2 \ar[rrr]^(.3){f_2}                       &&& \mathbb{C}^4/G \\
                     &&& Y_3 \ar[rrru]_{f_3}               && &
}
\]

\subsection{Singular}
Let $G$ is follwing type, then $\mathbb{C}^4/G$ has not crepant resolution, and ${\rm Hilb}^G(\mathbb{C}^4)$ has terminal singularity. This is a high-dimensional specific example.
\begin{Prop}\label{main3}
If $G$ is generated by $g=\frac{1}{2m}(1,2m-1,m,m)$, then ${\rm Hilb}^G(\mathbb{C}^4)$ has singularity.  
\end{Prop}

{\it Proof. }
To prove this proposition, there are two steps. The first step is to show that $\Gamma(I)$ is $G$-graph, where $\Gamma(I)$ is defined by $I=(x^m,y^m,xy,z^2,xz,yz,w)$. Next, the affine toric variety determined by $\sigma(\Gamma(I))$ has singulality.\\
We define the group homomorphism $\rho:G \to GL(1,\mathbb{C})$ such that it holds $\rho(g)=\varepsilon_{2m}^i$ for $g=\frac{1}{r}(i,2m-i,mi,mi) \in G$. We can see that $\rho$ is irreducible representation of $G$.
Let ${\rm Irr}(G) $ be the set of irreducible representations of $G$, then it is generated by $\rho$. In other word, we can write ${\rm Irr}(G)=\{\rho^0, \rho, \rho^2, \dots, \rho^{2m-1}\}$.\\
There is the bijection to ${\rm Irr}(G)$ from $\Gamma(I)$.
\begin{center}
\begin{eqnarray}
\Gamma(I) &\longrightarrow& {\rm Irr}(G) \nonumber \\
1  &\longmapsto& p^0 \nonumber \\
x &\longmapsto&  p  \nonumber \\
x^2 &\longmapsto&  p^2  \nonumber \\ 
     &\vdots&  \nonumber \\
x^{m-1} &\longmapsto&  p^{m-1}  \nonumber \\
z &\longmapsto&  p^m  \nonumber \\
y^{m-1} &\longmapsto&  p^{m+1}  \nonumber \\
y^{m-2} &\longmapsto&  p^{m+2}  \nonumber \\
         &\vdots&   \nonumber \\
y &\longmapsto&  p^{2m-1}  \nonumber 
\end{eqnarray}
\end{center}
Thus, $\Gamma(I)$ is $G$-graph. In addition, the following equation holds.
\begin{eqnarray}
{\rm wt}(x^m)&=&{\rm wt}(y^m)=p^m={\rm wt}(z) \nonumber \\
{\rm wt}(xz)&=&p^{m+1}={\rm wt}(y^{m-1}) \nonumber \\
{\rm wt}(yz)&=&p^{m-1}={\rm wt}(x^{m-1}) \nonumber \\
{\rm wt}(z)&=&{\rm wt}(w) \nonumber
\end{eqnarray}
From the above equation and $I=(x^m,y^m,xy,z^2,xz,yz,w)$, 
the coordinate ring of $\sigma(\Gamma(I))$ is $\mathbb{C}\left[\frac{x^m}{z},\frac{y^m}{z},\frac{xz}{y^{m-1}},\frac{yz}{x^{m-1}},\frac{w}{z}\right]\cong \mathbb{C}[X,Y,Z,W,V]/(XW-YZ)$. The affine variety $\textbf{V}(XW-YZ) \subseteq \mathbb{C}^5$ has terminal singularity. This means that the affine toric variety determined by $\Gamma(I)$ is singular.
\qed

\textsc{Graduate School Of Mathematical Sciences, University Of Tokyo 3-8-1 Komaba, Meguro-ku, Tokyo 153-8914, Japan}\\
E-mail:yusuke.sato@ipmu.jp
\end{document}